\documentclass{amsart}
\usepackage{amssymb,amscd,amsthm,amsmath}
\usepackage{comment}
\usepackage{epsfig}
\usepackage{graphicx}
\usepackage{xypic}
\usepackage{url}
\usepackage{mathtools}
\usepackage{hyperref}
\usepackage{quiver}
\usepackage{mathrsfs}

\numberwithin{equation}{section}

\makeatletter
\newcommand{\oset}[3][0ex]{%
 \mathrel{\mathop{#3}\limits^{
 \vbox to#1{\kern-2\ex@
 \hbox{\(\scriptstyle#2\)}\vss}}}}
\makeatother

\newcommand{\PP}{\mathbb{P}}
\newcommand{\OO}{\mathcal{O}}

\newcommand{\avgop}{\operatorname{Av}}
\newcommand{\bfe}{{\bf e}}

\newcommand{\Hom}{\textsf{Hom}}

\newcommand{\Spec}{\operatorname{Spec}}

\newcommand{\End}{\operatorname{End}}

\newcommand{\rk}{\operatorname{rk}}

\newcommand{\Ext}{\operatorname{Ext}}
\newcommand{\defe}{\operatorname{def}}

\def\Spec{\mathop{\mathrm{Spec}}}
\def\lSpec{\mathop{\mathcal Spec}\nolimits}
\newcommand{\Proj}{\operatorname{Proj}}

\newcommand{\SL}{\operatorname{SL}}

\newcommand{\defi}[1]{\textsf{#1}} 

\newtheorem{theorem}{Theorem}[section]
\newtheorem{lemma}[theorem]{Lemma}
\newtheorem{proposition}[theorem]{Proposition}
\newtheorem{corollary}[theorem]{Corollary}

\theoremstyle{definition}
\newtheorem{definition}[theorem]{Definition}

\newtheorem{remark}[theorem]{Remark}

\newtheorem{problem}[theorem]{Problem}
\newtheorem{example}[theorem]{Example}
\newtheorem{notation}[theorem]{Notation}

\definecolor{ao(english)}{rgb}{0.0, 0.5, 0.0}

\input xy
\xyoption{all}

\makeatletter
\def\dasharrowfill@#1#2#3#4{%
 $\m@th
 \thickmuskip0mu
 \medmuskip\thickmuskip
 \thinmuskip\thickmuskip
 \relax
 #4#1\mkern2mu
 \xleaders\hbox{$#4\mkern2mu#2\mkern2mu$}\hfill
 \mkern2mu
 #3$%
}

\def\dashrightarrowfill@{\dasharrowfill@\relbar\relbar\rightarrow}

\providecommand*\xdashrightarrow[2][]{%
  \ext@arrow 0055{\dashrightarrowfill@}{#1}{#2}}
\makeatother

\title{Ulrich bundles exist for sufficiently large multiples of any polarization}

\author{Izzet Coskun}
\address{Department of Mathematics, Statistics, and CS \\
University of Illinois Chicago, Chicago IL 60607}
\email{icoskun@uic.edu}
\author{Howard Nuer}
\address{Faculty of Mathematics \\
Technion, Israel Institute of Technology}
\email{hnuer@technion.ac.il}


\date{}

\begin{document}

\begin{abstract}
We prove that an integral complex projective variety $X$ admits an Ulrich bundle in some projective embedding if and only if $X$ is Cohen--Macaulay. In the Cohen--Macaulay case we answer the asymptotic existence question of Eisenbud and Schreyer by constructing Ulrich bundles for $mH$ for every ample Cartier divisor $H$ and every sufficiently large integer $m$. Examples of Anghel and Lopez show why one must allow the polarization to vary. Dropping either reducedness or irreducibility can produce a Cohen--Macaulay projective scheme with an ample divisor none of whose positive multiples admits an Ulrich bundle. For sheaves, we obtain asymptotic existence with full support on every reduced complex projective surface of pure dimension two. Anghel's obstruction to small Cohen--Macaulay modules gives normal projective threefolds with no full-support Ulrich sheaf in any embedding. If $X$ is integral and Cohen--Macaulay of dimension at least two, then $(X,mH)$ is Ulrich-wild for every sufficiently large $m$. For each such $m$, there is a constant $c>0$ and, for unbounded ranks $r$, a family of dimension at least $cr^2$ of pairwise non-isomorphic indecomposable Ulrich bundles of rank $r$ with the same determinant. The proof uses finite flat covers and noncommutative rank to construct a bundle map inducing simultaneous isomorphisms on cohomology. Multiplication maps give Ulrich bundles for every ample polarization on abelian varieties and smooth projective toric varieties, and Ulrich sheaves on singular normal projective toric varieties.
\end{abstract}

\maketitle

\section{Introduction}

Eisenbud and Schreyer asked whether every embedded projective variety is the support of an Ulrich sheaf \cite[p.~543]{ES03}. Examples of Anghel and Lopez show that even a smooth complex projective variety can have a very ample polarization for which it admits no Ulrich bundle \cite{Anghel,Lopez}. We consider how the answer changes when one takes multiples of an ample divisor.

\begin{definition}
Let $X$ be a complex projective scheme of pure dimension $n$ and let $H$ be an ample Cartier divisor on $X$. A nonzero coherent sheaf $E$ of dimension $n$ is an \defi{Ulrich sheaf} for $H$ if
\[
 H^j(X,E(-iH))=0 \qquad (1\le i\le n,\ 0\le j\le n).
\]
We say that $E$ has \defi{full support}, or \defi{support $X$}, if its scheme-theoretic support is $X$, that is, if $\operatorname{Ann}_{\OO_X}(E)=0$. When $X$ is reduced, this means that $E$ is nonzero at the generic point of every irreducible component. An \defi{Ulrich bundle} is a locally free Ulrich sheaf of positive rank.
\end{definition}

Ulrich bundles connect the geometry of projective varieties with
commutative algebra and explicit equations. When $H$ is very ample and $n>0$,
the graded module of twisted sections of an Ulrich bundle is a maximally
generated maximal Cohen--Macaulay module over the homogeneous coordinate ring
of $X$, with a linear resolution over the ambient polynomial ring
\cite[Proposition~2.1(d), p.~548]{ES03}. For a graded maximal Cohen--Macaulay module generated in degree zero over a
standard graded Cohen--Macaulay algebra, such a linear resolution is equivalent
to equality of the minimal number of generators and the multiplicity
\cite[Proposition~1.5, pp.~184--185]{BHU}.

A rank-$r$ Ulrich bundle on an embedded variety gives a determinantal formula for the $r$th power of its Chow form \cite[Theorem~0.3, p.~542]{ES03}. On a smooth cubic surface, the possible first Chern classes of Ulrich bundles are related to the minimal resolution conjecture for points on curves on the surface \cite[Corollary~1.4]{CKM}.
In Boij--S\"{o}derberg theory, Eisenbud and Schreyer describe the cohomology tables of vector bundles on projective space in terms of supernatural bundles \cite[Theorem~0.5]{ES09}. Write $C(X,L)$ for the cone of nonnegative rational combinations of the tables $(h^i(X,\mathcal F(jL)))_{i,j}$ of coherent sheaves $\mathcal F$. For a projective $n$-fold and a very ample divisor $L$, this cone equals that of $\PP^n$ if and only if $(X,L)$ admits an Ulrich sheaf of dimension $n$ \cite[Theorem~4.2, p.~43]{ES11}.

There is a large literature constructing Ulrich bundles on specific varieties. For example, existence is known for smooth curves \cite[Theorem~4.3 and Proposition~4.4, pp.~558--559]{ES03}, smooth projective complete intersections with their hyperplane class \cite{BaHU} (see \cite[Theorem~2.7]{Beauville}), Grassmannians in the Pl\"ucker embedding \cite{CostaMiroRoig1} and some two-step flag varieties with respect to the minimal very ample class \cite{CCHMW}.
Among rational surfaces, examples include Hirzebruch surfaces \cite{Antonelli,CH}, anticanonically embedded del Pezzo surfaces \cite[Corollary~6.5]{ES03} and certain blowups of the plane \cite[Theorem~6.6]{ES03}. Kim treats rational surfaces with an anticanonical pencil under a condition on the Clifford index \cite{Kim}.
For a smooth projective toric surface, Casnati's theorem gives an Ulrich bundle of rank two for every ample polarization \cite[Theorem~1.4]{CasnatiSurfaces}. Indeed, such a surface is rational, its toric boundary is anticanonical, and every ample line bundle on it is very ample.
Rank-two Ulrich bundles exist for every very ample polarization on a K3 surface \cite[Theorem~1]{Faenzi2}, an abelian surface \cite[Theorem~1]{BeauvilleU} or an Enriques surface \cite[Proposition~5.1]{Beauville}. For further constructions and results on Chern classes and stability on these surfaces, see \cite{AGO,Casnati,CNY} for K3 surfaces, \cite[Theorem~5.2]{CNY2} for abelian surfaces and \cite[Theorem~3.2]{BorisovNuer} for unnodal Enriques surfaces with a Fano polarization. We refer the reader to \cite{ecoskun} for a survey of higher-dimensional examples and basic properties of Ulrich bundles.

Our main theorem partially answers the original question of Eisenbud and Schreyer for arbitrary $X$ under some mild assumptions.

\begin{theorem}\label{thm-main}
Let $X$ be an integral Cohen--Macaulay complex projective variety and let $H$ be an ample Cartier divisor on $X$. Then $X$ admits an Ulrich bundle for $mH$ for every sufficiently large integer $m$.
\end{theorem}

Coskun and Huizenga proved an analogue of Theorem \ref{thm-main} for smooth surfaces \cite[Theorem~4.3, p.~43]{CH}: for a smooth polarized surface $(X,H)$, there exist Ulrich bundles of rank 2 on $X$ for $mH$ provided that $m$ is a sufficiently large positive integer. Their proof uses O'Grady's irreducibility theorem for moduli spaces of sheaves of large discriminant on a smooth projective surface \cite[Theorem~D]{ogrady}.  As such tools are unavailable in higher dimension, we prove Theorem \ref{thm-main} via constructive methods that we describe at the end of this section.  

\begin{remark}\label{rem-anghel-threshold}
The threshold in Theorem~\ref{thm-main} depends on the polarized variety, even for smooth surfaces. For a smooth complex surface $S$ with a very ample divisor $H$, Anghel's necessary condition for an Ulrich bundle for $mH$ is $m^2H^2\ge\sigma(S)$, where $\sigma(S)$ is the signature \cite[Proposition~2.1]{Anghel}. For his Hesse pairs $(Y_n,H_n)$, $n\ge3$, the ratio $\sigma(Y_n)/H_n^2=(3n^2-11)/7$ is unbounded \cite[Theorem~1.2]{Anghel}. Thus, for any $M\ge 1$, one can choose a single surface for which none of $H,2H,\ldots,MH$ admits an Ulrich bundle. In fact, whenever $7m^2<3n^2-11$, $(Y_n,mH_n)$ carries no nonzero arithmetically Cohen--Macaulay bundle, and its section ring admits no nonzero finitely generated graded maximal Cohen--Macaulay module \cite[Theorem~1.1 and Remark~4.2]{AnghelACM}. Lopez uses Anghel's surfaces to construct smooth examples in every dimension $n\ge 3$ for which a given very ample polarization admits no Ulrich bundle \cite[Theorem~1]{Lopez}.
\end{remark}

Normalizing the irreducible components gives an asymptotic answer to the sheaf question on every reduced projective surface of pure dimension two.

\begin{corollary}\label{cor-surface-sheaves}
Let $S$ be a reduced complex projective scheme of pure dimension two and let $H$ be an ample Cartier divisor on $S$. For every sufficiently large integer $m$, there is an Ulrich sheaf for $mH$ with support $S$.
\end{corollary}

\begin{proof}
Let $\nu_\alpha:\widetilde S_\alpha\to S$ be the finite map from the normalization of each irreducible component $S_\alpha$. The normal surface $\widetilde S_\alpha$ is Cohen--Macaulay by Serre's condition $S_2$, and $\nu_\alpha^*H$ is ample. By Theorem~\ref{thm-main}, there is an integer $m_\alpha$ such that $\widetilde S_\alpha$ admits an Ulrich bundle $E_{\alpha,m}$ for $m\nu_\alpha^*H$ whenever $m\ge m_\alpha$. For $m\ge\max_\alpha m_\alpha$, set
\[
 \mathcal F_m=\bigoplus_\alpha(\nu_\alpha)_*E_{\alpha,m}.
\]
Each summand has annihilator the ideal of $S_\alpha$, since $E_{\alpha,m}$ has positive rank and $\nu_\alpha$ is dominant onto $S_\alpha$. The intersection of these ideals is zero because $S$ is reduced. Thus $\mathcal F_m$ has full support. Finiteness and the projection formula give
\[
 H^j(S,\mathcal F_m(-imH))
 \cong\bigoplus_\alpha H^j(\widetilde S_\alpha,E_{\alpha,m}(-im\nu_\alpha^*H))=0
 \qquad(i=1,2,\ j\ge 0).\qedhere
\]
\end{proof}

Anghel asks how the cone of cohomology tables changes as one takes multiples of a polarization \cite[\S12, question~4]{AnghelSurvey}. Theorem~\ref{thm-main} gives
\[
 C(X,mH)=C(\PP^n,\OO_{\PP^n}(1))\qquad (m\gg0)
\]
for every integral Cohen--Macaulay complex projective $n$-fold $X$ and ample Cartier divisor $H$. The equality also holds on every reduced complex projective surface of pure dimension two by Corollary~\ref{cor-surface-sheaves}. Indeed, for a finite linear projection $\pi:X\to\PP^n$ in the embedding by $mH$, pushforward preserves cohomology tables, and an Ulrich sheaf $U$ has $\pi_*U\simeq\OO_{\PP^n}^{\oplus a}$ for some $a>0$. For every coherent sheaf $G$ on $\PP^n$, the table of $U\otimes\pi^*G$ is therefore $a$ times the table of $G$. When $X$ is Cohen--Macaulay and $U$ is a bundle, $\pi$ is flat and the same argument gives equality of the cones of vector bundle tables.

For bundles, Theorem~\ref{thm-main} gives a characterization in every dimension.
\begin{corollary}\label{cor-embedding}
For an integral complex projective variety $X$, the following are equivalent.
\begin{enumerate}
\item $X$ is Cohen--Macaulay.
\item $X$ admits an Ulrich bundle in some projective embedding.
\item For every ample Cartier divisor $H$ on $X$, there is an Ulrich bundle for $mH$ for every sufficiently large integer $m$.
\end{enumerate}
\end{corollary}

\begin{proof}
Theorem~\ref{thm-main} gives $(1)\Rightarrow(3)$, and $(3)\Rightarrow(2)$ since $mH$ is very ample for $m$ sufficiently large. For $(2)\Rightarrow(1)$, an Ulrich sheaf for a very ample polarization is locally Cohen--Macaulay on its support \cite[Proposition~2.1(d), p.~548]{ES03}. At each point $x$, a free $\OO_{X,x}$-module of positive rank has the same depth as $\OO_{X,x}$.
\end{proof}
The depth argument above makes Cohen--Macaulayness necessary at every point of a projective scheme carrying a positive-rank Ulrich bundle in an embedding. Example~\ref{ex-abelian-cone} illustrates this obstruction. Examples~\ref{ex-square-zero-thickening} and~\ref{ex-four-planes} show that dropping reducedness or irreducibility, respectively, can give a Cohen--Macaulay scheme with an ample divisor none of whose positive multiples admits an Ulrich bundle.
\begin{example}[Cohen--Macaulayness]\label{ex-abelian-cone}
Theorem~\ref{thm-main} applies in particular to integral projective local complete intersections and to normal projective surfaces. In higher dimension, however, normal projective varieties that fail to be Cohen--Macaulay cannot admit an Ulrich bundle in any projective embedding. The projective cone over a projectively normal abelian surface $A$ is a normal threefold. Since $H^1(A,\OO_A)\ne 0$, the local ring at the vertex has depth two, so Corollary~\ref{cor-embedding} rules out an Ulrich bundle in every projective embedding. Indeed, for the homogeneous coordinate ring $R$ of $A$ and its irrelevant maximal ideal $\mathfrak m$, the degree-zero part of $H^2_{\mathfrak m}(R)$ is $H^1(A,\OO_A)$, while normality gives $\operatorname{depth} R_{\mathfrak m}\ge 2$.

The cone carries an Ulrich sheaf with full support in its cone embedding. Take Beauville's rank-two Ulrich bundle $U$ on $A$ \cite[Theorem~1]{BeauvilleU} and let $M=\bigoplus_{j\in\mathbb Z}H^0(A,U(j))$. The graded $R[t]$-module $M[t]$ defines the desired sheaf on $\Proj R[t]$. Since $U$ has positive rank, $M[t]$ has positive rank over the domain $R[t]$ and hence full support. Adjoining $t$ preserves Cohen--Macaulayness and extends the linear resolution of $M$ over the ambient polynomial ring \cite[Proposition~2.1(d)]{ES03}.
\end{example}

Corollary~\ref{cor-surface-sheaves} fails for normal threefolds. Let $R$ be the section ring of a Hesse pair $(Y_n,H_n)$ with $n\ge3$. Anghel proves that the completion of $R$ at the irrelevant maximal ideal $R_+$ admits no small Cohen--Macaulay module, that is, no nonzero finitely generated maximal Cohen--Macaulay module \cite[Corollary~9.3]{AnghelSurvey}. Since $Y_n$ is smooth, its section ring is normal, so $X=\Proj R[t]$, with $\deg t=1$, is a normal projective threefold. The chart $D_+(t)$ is $\Spec R$, with vertex $v$ corresponding to $R_+$. In any projective embedding of $X$, a full-support Ulrich sheaf would be locally Cohen--Macaulay by \cite[Proposition~2.1]{ES03}. Its completed stalk at $v$ would therefore be a nonzero finite module of depth three over $\widehat{\OO}_{X,v}=\widehat{R_{R_+}}$, contradicting Anghel's result.

\begin{example}[Reducedness]\label{ex-square-zero-thickening}
Let $p:X\to\PP^n$, $n\ge 1$, be the split square-zero thickening with algebra $\OO_{\PP^n}\oplus\OO_{\PP^n}(1)$, and set $\OO_X(H)=p^*\OO_{\PP^n}(1)$. The map $p$ is finite flat, $H$ is ample, and $X$ is locally the hypersurface $z^2=0$, hence Gorenstein. Write $i:\PP^n=X_{\mathrm{red}}\hookrightarrow X$ for the closed immersion. The sheaf $i_*\OO_{\PP^n}$ is Ulrich for $H$ since $p\circ i$ is the identity.

In fact every $n$-dimensional Ulrich sheaf $E$ for any $mH$, $m\ge1$, has scheme-theoretic support $X_{\mathrm{red}}$. Finiteness and the projection formula make $V=p_*E$ Ulrich for $\OO_{\PP^n}(m)$. Since $\PP^n$ is smooth, $V$ is locally free by the local Cohen--Macaulay property, and it is slope-semistable by \cite[Theorem~2.9(a)]{CHGS}. The action of the square-zero ideal is a map $V(1)\to V$. Both bundles are semistable and $\mu(V(1))>\mu(V)$, so the map vanishes. Hence $E=i_*V$, whose annihilator is exactly the square-zero ideal. In particular no $mH$ admits an Ulrich bundle on $X$. For comparison, on the split double structure with algebra $\OO_{\PP^n}[\varepsilon]/(\varepsilon^2)$, the structure sheaf is Ulrich for the pullback of $\OO_{\PP^n}(1)$, since its pushforward is $\OO_{\PP^n}^{\oplus 2}$.
\end{example}

\begin{example}[Irreducibility]\label{ex-four-planes}
Irreducibility is needed in Theorem \ref{thm-main} even for a connected, reduced Cohen--Macaulay surface. Let $S_0=\PP^2$. The map $g:S_0\to\PP^2$ given by $[x:y:z]\mapsto[x^2:xy+z^2:y^2]$ is finite of degree four, since the quadrics have no common zero and $g^*\OO_{\PP^2}(1)=\OO_{S_0}(2)$. It maps the line $z=0$ isomorphically onto a smooth conic $D$. Glue $S_0$ to three copies $S_1,S_2,S_3$ of the target by identifying this line with $D$ in each copy. The resulting surface has a finite map $p:X\to\PP^2$ with $X=\lSpec_{\PP^2}\mathcal A$, where $\mathcal A$ is the subalgebra of $g_*\OO_{S_0}\times\OO_{\PP^2}^3$ consisting of tuples whose restrictions to $D$ agree. Thus $X$ is connected and reduced. Projection onto the first entry gives a surjection $\mathcal A\to g_*\OO_{S_0}$ with kernel $\OO_{\PP^2}(-2)^{\oplus3}$. Since $g$ is flat by miracle flatness \cite[Tag 00R4]{Stacks}, $\mathcal A$ is locally free of rank seven. Hence $p$ is finite flat over the smooth plane, $X$ is Cohen--Macaulay and projective, and $\OO_X(H)=p^*\OO_{\PP^2}(1)$ is ample.

Let $E$ be a vector bundle of rank $r>0$ on $X$, and write $E_i=E|_{S_i}$. Since the determinants agree on $D$, their degrees on the four planes are $2a,a,a,a$ for some integer $a$. Restriction to $S_0$ gives an exact sequence on the target
\[
0\longrightarrow\bigoplus_{i=1}^3 E_i(-2)\longrightarrow p_*E\longrightarrow g_*E_0\longrightarrow0.
\]
Sections on $S_0$ vanishing along the gluing line extend by zero and give a subsheaf $g_*(E_0(-1))\subset p_*E$. For a rank-$r$ bundle $V$ on $S_0$, comparing linear coefficients in $\chi((g_*V)(t))=\chi(V(2t))$ by Riemann--Roch gives $\deg c_1(g_*V)=2\deg c_1(V)-3r$. The exact sequence therefore gives
\[
\mu\bigl(g_*(E_0(-1))\bigr)-\mu(p_*E)
=\frac{4a-5r}{4r}-\frac{7a-9r}{7r}=\frac1{28}
\]
with slopes measured against $\OO_{\PP^2}(1)$. This subsheaf destabilizes $p_*E$ for every $\OO_{\PP^2}(m)$ with $m\ge1$. If $E$ were Ulrich for $mH$, the projection formula would make $p_*E$ Ulrich for $\OO_{\PP^2}(m)$ and hence slope-semistable \cite[Theorem~2.9(a)]{CHGS}. Thus no positive multiple of $H$ admits an Ulrich bundle.

For every $m\ge1$, the restrictions of $mH$ to the four planes are $\OO_{S_0}(2m)$ and $\OO_{S_i}(m)$ for $1\le i\le3$. Each polarized plane carries an Ulrich bundle \cite[Corollary~5.7]{ES03}. Pushing their direct sum forward from the normalization $\bigsqcup_{i=0}^3 S_i$ gives an Ulrich sheaf with full support for $mH$, by the proof of Corollary~\ref{cor-surface-sheaves}.
\end{example}

For special $X$, we can sometimes use the asymptotic theorem to give an Ulrich bundle for the original
polarization $L$. Suppose that $X$ is as in Theorem~\ref{thm-main} and that
$f:X\to X$ is finite flat with $f^*L\simeq L^{\otimes q}$ for an ample
line bundle $L$ and an integer $q>1$. The $s$th iterate $f^s$ pulls
$L$ back to $L^{\otimes q^s}$. For $s$ sufficiently large,
Theorem~\ref{thm-main} gives an Ulrich bundle $U$ for this power,
and Lemma~\ref{lem-passtocover} shows that $(f^s)_*U$ is Ulrich for $L$.
Its rank is $(\deg f)^s\rk U$, where $\deg f=q^{\dim X}$ by the
intersection formula for the top power of $L$.

\begin{corollary}\label{cor-multiplication}
Let $X$ be a complex abelian variety or a smooth complex projective
toric variety. Then $X$ admits an Ulrich bundle for every ample line
bundle $L$ on $X$.
\end{corollary}

\begin{proof}
First let $X$ be an abelian variety of dimension $g$. For $m\geq2$,
multiplication by $m$ is finite \'etale of degree $m^{2g}$. Set
$M=L^{\otimes m^2}$. If $L$ is symmetric, $[m]^*L\simeq M$. For
arbitrary $L$, the two line bundles have the same class in the
N\'eron--Severi group \cite[Corollary~6.6 and Theorem~8.2]{MilneAV}.
Since $M$ is ample, the map $x\mapsto t_x^*M\otimes M^{-1}$ from
$X$ to $\operatorname{Pic}^0(X)$ is surjective
\cite[Proposition~10.1]{MilneAV}. Choose $x$ such that
$t_x^*M\simeq[m]^*L$. Then $f=[m]\circ t_{-x}$ is finite \'etale
and satisfies $f^*L\simeq M$, so the preceding argument applies.

Now let $X$ be toric. The $m$th power map on the dense torus extends
to a finite surjective endomorphism $F_m:X\to X$
\cite[\S 2]{FujinoToric}. Represent $L$ by a torus-invariant Cartier
divisor. Its local equations are characters, which $F_m$ pulls back
to their $m$th powers. Hence $F_m^*L\simeq L^{\otimes m}$. On the
affine chart of any smooth maximal cone, $F_m$ sends each coordinate
$x_i$ to $x_i^m$. The monomials $\prod_{i=1}^{n}x_i^{a_i}$ with
$0\leq a_i<m$ form a free basis over $\mathbb C[x_1^m,\ldots,x_n^m]$,
where $n=\dim X$. These charts cover $X$ when $X$ is smooth, so
$F_m$ is finite flat and the result follows.
\end{proof}

\begin{remark}
This argument also gives an Ulrich sheaf for every ample line bundle $L$ on a singular normal projective toric variety $X$. Indeed, $X$ is Cohen--Macaulay \cite[Theorem~9.2.9]{CLS}. So we may take an Ulrich bundle $U$
for $L^{\otimes m}$ with $m$ sufficiently large by Theorem \ref{thm-main}. Finiteness of $F_m$
and the projection formula give
\[
 H^j(X,(F_m)_*U\otimes L^{-i})\cong H^j(X,U\otimes L^{-mi})=0
 \qquad (1\leq i\leq\dim X,\; j\geq0).
\]
Thus $(F_m)_*U$ is an Ulrich sheaf with support $X$. At a singular
point $F_m$ can fail to be flat, so the pushforward need not be
locally free.  
\end{remark}

We do not attempt to give an effective bound for the rank of the Ulrich bundle constructed in Theorem \ref{thm-main} or for the multiple of the polarization needed. This raises the following two problems.

\begin{problem}
Let $(X,H)$ be a polarized, smooth, complex projective variety.
    \begin{enumerate}
    \item Determine bounds on the multiple $m$ of $H$ needed to guarantee the existence of an Ulrich bundle for $(X, mH)$.
     \item Determine bounds on the ranks of the Ulrich bundles that exist for $(X,mH)$.
    \end{enumerate}
    Of course, constructing Ulrich bundles on special varieties for specific polarizations remains an interesting problem.
\end{problem}

For part~(1), Remark~\ref{rem-anghel-threshold} gives the lower bound $m^2H^2\ge\sigma(X)$ on surfaces when $mH$ is very ample. More generally, if $n=\dim X\ge2$ and $mH$ is very ample, Lopez's inequality for an Ulrich bundle on $(X,mH)$ \cite[Lemma~3.1]{Lopez} reads
\[
 (n+1)m^2H^n\ \ge\ (K_X^2-2c_2(X))\cdot H^{n-2}.
\]
Anghel proves the same bound for any nonzero arithmetically Cohen--Macaulay bundle when $mH$ is ample and base-point-free \cite[Proposition~5.1]{AnghelACM}. For a fixed Hesse pair, Anghel asks which multiples between the excluded range and the range of asymptotic existence admit arithmetically Cohen--Macaulay bundles \cite[Remark~6.1]{AnghelACM}. The analogous question for Ulrich bundles is an instance of part~(1). Since every two-dimensional Ulrich sheaf on the smooth surface $Y_n$ is locally free, $C(Y_n,mH_n)$ is strictly smaller than $C(\PP^2,\OO_{\PP^2}(1))$ when $7m^2<3n^2-11$. The cones agree for all sufficiently large $m$, and determining which intermediate multiples give equality remains open.

Once existence is known, another well-studied question asks how large families of pairwise non-isomorphic Ulrich bundles can be.  More specifically, we have the following definition.
\begin{definition}
   A polarized projective scheme $(X,H)$ is of \defi{Ulrich-wild representation type}, or simply, \defi{Ulrich-wild}, if $X$ admits families of arbitrarily large dimension parameterizing pairwise non-isomorphic indecomposable Ulrich bundles.
\end{definition}

Costa, Mir\'o-Roig, and Pons-Llopis proved that a product of two or more positive-dimensional projective spaces is Ulrich-wild in its Segre embedding, with the sole exception of the quadric surface $\PP^1\times\PP^1$ \cite[Theorem~4.6]{CMP}.

Varying the map in our construction gives a lower bound on the dimension of a family in terms of the rank.

\begin{theorem}\label{thm-ulrichwild}
Let $X$ be an integral Cohen--Macaulay complex projective variety of dimension at least two and let $H$ be an ample Cartier divisor. For every sufficiently large integer $m$, put $P=mH$. There is a constant $c_P>0$ such that, for an unbounded sequence of ranks $r$, there is a family of dimension at least $c_Pr^2$ of pairwise non-isomorphic indecomposable Ulrich bundles of rank $r$ for $P$. All bundles in each family have the same determinant. In particular, $(X,P)$ is Ulrich-wild.
\end{theorem}

\begin{remark}
 For smooth surfaces, Coskun and Huizenga proved Ulrich-wildness for all sufficiently large multiples of the polarization \cite[Corollary~4.5, p.~45]{CH}.
\end{remark}

On an integral projective curve, start with a line bundle $B$ with $H^1(B)=0$ and successively subtract smooth points where a remaining section is nonzero. Each point lowers $h^0$ by one and preserves the vanishing of $H^1$. The resulting line bundle $A$ is acyclic, so $A(mH)$ is Ulrich for $mH$. The zero-dimensional case is immediate. For this reason, and to avoid trivialities, for the rest of the paper we assume that $\dim(X)=n\ge 2$.

\subsection*{The strategy of the proof}
Our starting point is a finite projection $X\to\PP^n$, where $n=\dim X$. We look for a kernel presentation using twists of bundles pulled back from projective space. At each of the $n$ twists in the Ulrich condition, we can arrange that the source and target have cohomology in a single degree and that the two groups have the same dimension. The difficulty is to find one surjection inducing an isomorphism at all $n$ twists. ChatGPT suggested drawing on Huizenga's use of the power map $u\mapsto u^d$ on $\PP^1$ to transfer Kronecker semistability for multiplication maps \cite[proofs of Lemma~2.2 and Proposition~3.1]{HuizengaSteiner}, and on the work of Coskun, Huizenga and Smith, who use strong Kronecker stability to prove maximal rank on global sections in their study of kernel bundles on $\PP^n$ \cite[Theorems~2.3 and~2.6]{CHS}.

To control the induced maps, we pass to a finite flat cyclic cover $f:Y\to X$ of degree $\ell$, equipped with a second finite projection $\pi:Y\to\PP^n$. For the line bundles used below, each summand of degree $d$ in the pushforward from $X$ gives $\ell$ summands of consecutive degrees $d,d-1,\ldots,d-\ell+1$ in the pushforward from $Y$. In the corresponding graded modules, multiplication realizes truncations of length $\ell$ of the evaluation filtrations on a general zero-dimensional linear section of $X$. We choose $\ell$ large enough that every proper nonzero graded submodule of each truncation has strictly larger mean degree than the whole module (Proposition~\ref{prop-multiplication}). This inequality is what we need to control the ranks of the maps on cohomology.

Write $m=qM+a$, where $L=qH$ is very ample and $0\le a<q$. On $Y$ we seek a sequence
\[
0\longrightarrow K\longrightarrow F^{\oplus 2N}
 \xrightarrow{\;\psi\;} F(\delta f^*L)^{\oplus N}\longrightarrow 0.
\]
We take $F=f^*B_0\otimes\pi^*G$, where $B_0$ is a line bundle on $X$ and $G$ is a sum of homogeneous bundles on $\PP^n$. We choose $B_0$ so that $B_i=B_0(-iaH)$ has vanishing intermediate cohomology at every twist by $L$ for $0\le i\le n$. By Horrocks' criterion \cite{Horrocks}, each $\pi_*f^*B_i$ then splits into line bundles (Lemmas~\ref{lem-choosingta} and~\ref{lem-directsumoflb}). This splitting and the Borel--Bott--Weil theorem allow us to choose $G$ so that at the twist $-imf^*H$ both terms of $\psi$ have cohomology of the same dimension, concentrated in degree $i-1$ (Theorem~\ref{thm-constructbundle} and Proposition~\ref{prop-surjhom}).

We degenerate the cohomology maps to maps given by multiplication on the truncations. At each twist, the mean-degree inequality implies that the images of any subspace of the cohomology of one copy of $F$ under the limiting maps together span a space of at least twice its dimension (Proposition~\ref{prop-noshrunk}). The factor of two accounts for the $2N$ copies in the source and $N$ in the target. The noncommutative rank criterion of Fortin--Reutenauer \cite{FR} and the results of Derksen--Makam \cite{DM} then give invertible matrices of these maps for all sufficiently large $N$. Lifting through the degeneration gives a nonempty open subset of the space of bundle maps for each isomorphism condition. Everywhere surjective maps also form a nonempty open subset. A map $\psi$ in their intersection has Ulrich kernel $K$ by the long exact cohomology sequence. Since $f$ is finite flat, $f_*K$ is the desired Ulrich bundle on $X$ (Lemma~\ref{lem-passtocover}).

Given an Ulrich bundle for a very ample divisor $H_0$, Theorem~\ref{thm-wildrep} transfers families from projective space to suitable multiples of $H_0$. For wildness at every sufficiently large $m$, we set $A=f_*F$ and vary the map $A^{\oplus 2N}\to A(\delta L)^{\oplus N}$ on $X$. Lemma~\ref{lem-uniform-homs} lets us choose $\delta$ before taking $M$ large. With $m$ fixed and $N$ increasing, Lemma~\ref{lem-moving-kernels} gives a quadratic lower bound in $N$ for the dimension of the locus of framed kernels after subtracting the dimension of the group of frame changes. Keeping track of the ranks of their indecomposable summands in Lemma~\ref{lem-framed-ulrich} gives the quadratic bound of Theorem~\ref{thm-ulrichwild}. Fixing the determinant costs at most $\dim\operatorname{Pic}^0(X)$ parameters.

\begin{remark}
  We use characteristic zero at three points. In Lemma~\ref{lem-uniformposition} the monodromy group is the full symmetric group, in Proposition~\ref{prop-noshrunk} the relevant representations are completely reducible, and in \S\ref{sec-numerical} the Borel--Bott--Weil theorem computes the cohomology of homogeneous bundles. It would be interesting to construct Ulrich bundles in higher dimension for sufficiently large multiples of the polarization over an algebraically closed field of characteristic $p>0$.

  Over such a field, Pine constructs Ulrich bundles for every very ample polarization $H$ on a smooth projective surface $S$ \cite[Theorem~1.1]{Pine}. For each bundle $E$ in his construction, of rank $r$, he computes \cite[Remark~2.9]{Pine}
  \[
   \frac{2r\,c_2(E)-(r-1)c_1(E)^2}{r^2}
   =\frac{H^2+8\chi(\OO_S)-K_S^2}{4}.
  \]
  Thus these bundles violate the Bogomolov inequality precisely when $H^2<K_S^2-8\chi(\OO_S)$. On a smooth complex surface $K_S^2-8\chi(\OO_S)=\sigma(S)$, so this is exactly the numerical range in which Remark~\ref{rem-anghel-threshold} excludes Ulrich bundles.

  For a smooth projective $n$-fold over an algebraically closed field of arbitrary characteristic, Pine also constructs weakly Ulrich bundles of rank $n!$ for all sufficiently large multiples of a very ample divisor \cite[Theorem~1.3]{Pine}. For the original Hesse polarizations $(Y_n,H_n)$ over $\mathbb C$, Genc and Lopez prove that there are no weakly Ulrich bundles when $n\ge8$ \cite[Theorem~2]{GencLopez}.
\end{remark}

\subsection*{Acknowledgments} We thank Cristian Anghel for reading the preprint carefully and for his helpful comments, in particular for pointing out the obstruction to full-support Ulrich sheaves on the cones over the Hesse surfaces. We also thank Lev Borisov, Jack Huizenga and Matthew Woolf for invaluable discussions on Ulrich bundles. Izzet Coskun was partially supported by the NSF grant DMS-2200684 and the Simons Foundation Travel Support for Mathematicians award SFI-MPS-TSM-00013580; Howard Nuer was partially supported by ISF Grant 3175/25.

\subsection*{AI disclosure}
We proposed using a finite projection and developed the kernel construction. Two Astra models in OpenAI's ChatGPT suggested passing to the auxiliary cover $Y$. ChatGPT also suggested the normalization argument in Corollary~\ref{cor-surface-sheaves} and the fixed-determinant refinement of Theorem~\ref{thm-ulrichwild}, and helped edit the exposition. We take responsibility for the proofs, references and conclusions.

\section{Setup and preliminaries}\label{sec-prelim}
In this section, we introduce the setup of the proof.  We will preserve the notation introduced in this section throughout the paper. 
First, recall the following easy consequence of Serre vanishing and Serre duality.

\begin{lemma} \label{lem-easyconsserre}
    Let $X$ be an integral Cohen--Macaulay projective variety of dimension $n$ and let $H$ be an ample Cartier divisor. Then there exists an integer $C \geq 1$ such that
    $$H^j(X, \OO_X(mH)) = 0 \quad \mbox{for} \ 1 \leq j \leq n-1 \ \mbox{and}\ |m| \geq C.$$
\end{lemma}

\begin{proof}
    Since $H$ is ample, by Serre vanishing \cite[Ch.~III, Proposition~5.3, p.~229]{Hartshorne}, there exists $m_0$ such that $$H^j(X, \OO_X(mH)) = 0 \quad \mbox{for} \ j \geq 1 \ \mbox{if} \ m \geq m_0.$$
    By Serre duality for the Cohen--Macaulay variety $X$ with dualizing sheaf $\omega_X$ \cite[Ch.~III, Corollary~7.7, p.~244]{Hartshorne},
    $$H^j(X, \OO_X(mH))\cong H^{n-j}(X, \omega_X(-mH))^{\vee}.$$ By Serre vanishing, the latter group vanishes for $n-j \geq 1$ if $m \leq -m_1$ for some $m_1 > 0$.  If $C= \max(m_0, m_1)$, then
    $$H^j(X, \OO_X(mH)) = 0 \quad \mbox{for} \ 1 \leq j \leq n-1 \ \mbox{and}\ |m| \geq C.$$
\end{proof}

\begin{notation}
  Throughout the paper, let $p(t) = \chi(\OO_X(tH))$ be the Hilbert polynomial of $H$ on $X$. 
\end{notation}

\begin{remark}\label{rem-hpobs}
    For integers $\alpha$ and $q$ we will use the identity
    $$\chi(\OO_X((tq+\alpha)H))= p(qt+\alpha).$$
\end{remark}

\begin{notation}\label{def-C}
 Given $X$ and $H$, let $C$ be the integer determined by Lemma \ref{lem-easyconsserre}, so that $$H^j(X, \OO_X(mH))=0 \ \mbox{for} \ 1\leq j \leq n-1 \ \mbox{and} \ |m| \geq C.$$   
\end{notation}

\begin{definition}[Definition of $q$]\label{def-q} Fix an integer $q$ that satisfies the following three conditions:
\begin{enumerate}
    \item $q > (n+1)(2C-1)$,
    \item $qH$ is very ample on $X$, 
    \item $q$ is larger than the distance between any two roots of $p(t)$.
\end{enumerate}
Let $L=\OO_X(qH)$ and let $e = L^n$ be the degree of $X$ under the embedding by $L$. This notation will be fixed throughout the paper.
\end{definition}

\begin{lemma}\label{lem-choosingta}
    For each residue class $a \in \{ 0, 1, \dots, q-1\}$ modulo $q$, there exists an integer $t_a \in \{ 0, 1, \dots, q-1\}$ such that for any integer $k$, we have
$$H^j(X, \OO_X((t_a - ia+kq)H))=0 \quad \mbox{for} \ 1 \leq j \leq n-1 \ \mbox{and} \ 0 \leq i \leq n.$$ 
\end{lemma}
\begin{proof}
  Since $H^j(X, \OO_X(mH))=0$ for $1 \leq j \leq n-1$ and $|m| \geq C$, there are at most $2C-1$ residue classes mod $q$ that we need the twists $t_a - ia$ to avoid.  If we twist these $2C-1$ integers by $ia$ for $0 \leq i \leq n$, we get at most $(n+1)(2C-1)$ integers. Since $q> (n+1)(2C-1)$, there must exist an integer $0 \leq t_a \leq q-1$, which avoids them all. Twisting by integer multiples of $q$ does not change the residue class mod $q$, hence the lemma follows.
\end{proof}

\begin{notation}\label{not-LBi}
For the rest of the paper, for each residue class $a$ mod $q$, fix an integer $t_a$ satisfying the conclusions of  Lemma \ref{lem-choosingta}. Set $B_i = \OO_X((t_a -ia)H).$
\end{notation}

\begin{remark}\label{rem-fixedresidue}
    Write \(m=qM+a\), where \(0\leq a<q\). It suffices to prove the theorem separately in each residue class modulo \(q\). Indeed, suppose that for each \(a\) there is an integer \(M_a\) such that \(X\) admits an Ulrich bundle with respect to \(\mathcal O_X((qM+a)H)\) for every \(M\geq M_a\). If \(m\geq q\max_{0\leq b<q}M_b\), then \(M=\lfloor m/q\rfloor\geq M_a\) for the residue \(a\) of \(m\) modulo \(q\). Thus the theorem holds for every sufficiently large \(m\).

For the rest of the proof, fix an arbitrary \(a\in\{0,\ldots,q-1\}\) and write \(m=qM+a\). We will construct an Ulrich bundle for \(mH\) for every sufficiently large \(M\), allowing the construction and the lower bound on \(M\) to depend on \(a\).
\end{remark}

\begin{definition}\label{def-Y}
The cyclic cover associated with a nonzero section \(s\in H^0(X,L^\ell)\) is
\[
f\colon Y=\lSpec_X\!\left(\bigoplus_{k=0}^{\ell-1}L^{-k}\right)\longrightarrow X,
\]
where multiplication is induced by \(s\colon L^{-\ell}\to\mathcal O_X\). 
Then \(f\) is finite flat of degree \(\ell\). Set \(L_Y=f^*L\).
\end{definition}

In \S\ref{sec-cover} we choose \(s\) and construct a finite flat map \(\pi\colon Y\to\mathbb P^n\) with \(\pi^*\mathcal O_{\mathbb P^n}(1)\cong L_Y\).

\begin{lemma}\label{lem-passtocover}
    Let $f: Y \to X$ be a finite flat cover of $X$ of degree $\ell$. Let $A$ be an ample divisor on $X$ and let $D= f^* A$. If $V$ is an Ulrich bundle of rank $r$ on $Y$ for $D$, then $f_* V$ is an Ulrich bundle of rank $r \ell$ on $X$ for $A$.
\end{lemma}

\begin{proof}
    Since $f$ is finite flat of degree $\ell$, restriction of scalars takes $V$ to a vector bundle $f_*V$ of rank $r\ell$. Since $f$ is affine, its higher direct images vanish on quasi-coherent sheaves. Thus, by \cite[Ch.~III, Exercise~8.2, p.~252]{Hartshorne} or \cite[Tags 01XC and 089W]{Stacks},
\begin{equation}\label{eq-usef}
H^j (Y, V(-iD))= H^j(X, f_*(V(-iD)) ) \ \quad \mbox{for all} \ j.
\end{equation}
By the projection formula \cite[Ch.~III, Exercise~8.3, p.~253]{Hartshorne} or \cite[Tag 01E8]{Stacks},
\begin{equation}\label{eq-projectionf}
    f_*(V(-iD)) = f_*(V \otimes_{\OO_Y} f^*\OO_X(-iA)) \cong f_*V \otimes \OO_X(-iA) = (f_*V)(-iA).
\end{equation}
Since $V$ is an Ulrich bundle on $Y$ for $D$, we have $H^j (Y, V(-iD))=0$ for $1 \leq i \leq n$ and all $j$. Combining Equations \eqref{eq-usef} and \eqref{eq-projectionf}, we conclude that $$H^j(X, f_*V (-iA))=0 \ \mbox{for} \ 1 \leq i \leq n \ \mbox{and all} \ j.$$ Thus $f_*V$ is an Ulrich bundle on $X$ for $A$. 
\end{proof}

By Lemma \ref{lem-passtocover}, in order to construct an Ulrich bundle on $X$ for $\OO_X(mH)$, it suffices to construct an Ulrich bundle on $Y$ for $f^* \OO_X(mH)$. For the rest of the paper, we will work on $Y$.

\begin{definition}\label{def-Ai}
    Let $f: Y \to X$ be a cyclic cover as in Definition \ref{def-Y} and let $$A_i = f^* B_i = f^* \OO_X((t_a-ia)H).$$ 
\end{definition}

\begin{lemma}\label{lem-directsumoflb}
    Let $\pi: Y \to \PP^n$ be a finite flat map such that $\pi^* \OO_{\PP^n}(1) = L_Y$. Then $\pi_* A_i$ is a direct sum of line bundles for every $0 \leq i \leq n$.
\end{lemma}
\begin{proof}
    Since $\pi$ is finite flat and $A_i$ is a line bundle, \(\pi_*A_i\) is locally free. By Horrocks' splitting criterion \cite[Ch.~I, Thm.~2.3.1, p.~21]{OSS}, it suffices to show that \((\pi_*A_i)(k)\) has no intermediate cohomology for any \(k\in\mathbb Z\).

Both \(\pi\) and \(f\) are finite, so their higher direct images vanish on quasi-coherent sheaves. Applying the projection formula first to \(\pi\) and then to \(f\) gives
\[
\begin{aligned}
H^j\bigl(\PP^n,(\pi_*A_i)(k)\bigr)
&\cong H^j\bigl(Y,A_i\otimes L_Y^k\bigr)\\
&\cong \bigoplus_{h=0}^{\ell-1}
H^j\!\left(X,\OO_X\bigl((t_a-ia+(k-h)q)H\bigr)\right),
\end{aligned}
\]
where the second isomorphism uses \(A_i=f^*\OO_X((t_a-ia)H)\), \(L_Y=f^*\OO_X(qH)\), and \(f_*\OO_Y\cong\bigoplus_{h=0}^{\ell-1}\OO_X(-hqH)\).

Each integer \(t_a-ia+(k-h)q\) is congruent to \(t_a-ia\) modulo \(q\). By Lemma~\ref{lem-choosingta}, every summand therefore vanishes for \(1\le j\le n-1\), as required.
\end{proof}

\begin{notation}\label{not-directsum}
    By Lemma \ref{lem-directsumoflb}, $\pi_* A_i$ is a direct sum of line bundles, so set $$\pi_* A_i = \bigoplus_{u}\OO_{\PP^n}(-u)^{b_{iu}}.$$ We will adhere to this notation for the rest of the paper.
\end{notation}

Let $x_0, \dots, x_n \in H^0(X, L)$ be general sections of $L=\OO_X(qH)$. These define a finite map $\phi:X \to \PP^n$ of degree $e= L^n$. Since $X$ is Cohen--Macaulay and $\PP^n$ is regular of the same dimension, $\phi$ is flat by miracle flatness \cite[Tag 00R4]{Stacks}.

\begin{corollary}\label{cor-projBi}
   For each $0 \leq i \leq n$, there exist nonnegative integers $d_{iu}$, $1 \leq u \leq e$, such that the integers $d_{iu}$ are not all equal to each other and $$\phi_* B_i \cong \bigoplus_{u=1}^e \OO_{\PP^n}(-d_{iu}).$$ 
\end{corollary}
\begin{proof}
   By applying Lemma \ref{lem-directsumoflb} when the cyclic cover is trivial, we conclude that $\phi_* B_i$ is a direct sum of line bundles. Since $B_i(-qH)$ is a negative multiple of $\OO_X(H)$, it cannot have any global sections. We conclude that $H^0(\PP^n,\phi_*B_i(-1))=0$, hence all the summands occurring in $\phi_* B_i$ must be non-positive, yielding $d_{iu} \geq 0$. Finally, if all the integers $d_{iu}$ are equal, we conclude that the Hilbert polynomial of $\phi_* B_i$ is $p(t_a-ia+qt)$ and has roots which are $n$ consecutive integers. Consequently, the polynomial $p$ has roots that are $q$ apart, contrary to Condition (3) in Definition \ref{def-q} defining $q$. We conclude that at least two of the integers $d_{iu}$ must be distinct. 
\end{proof}

\begin{remark}
    The fact that the integers $d_{iu}$ in Corollary \ref{cor-projBi} are not all equal will be important in the next section in proving the slope inequality in Proposition \ref{prop-multiplication}.
\end{remark}

\section{Construction of the cover $Y$}\label{sec-cover}
In this section, we explain how to construct a cyclic cover $Y$ of $X$ so that multiplication maps satisfy certain rank inequalities. These will be important for proving cohomology vanishing later in the proof.
We preserve the notation from \S \ref{sec-prelim}.

\begin{notation}
Set $B$ to be one of  $B_i= \OO_X((t_a-ia)H)$ for $0 \leq i \leq n$. Choose general sections $x_0, \dots, x_n \in H^0(X, L)$ and let $\phi: X \to \PP^n$ be the finite projection of degree $L^n=e$ they define. Since $L$ is very ample and the sections are chosen generally, by Bertini's Theorem, the common zero locus of any $n$ of these sections is a reduced set of $e$ points. Let $Z := V(x_1, \dots, x_n)$ be the vanishing locus of the last $n$ sections.
\end{notation}

We recall that, by Corollary \ref{cor-projBi}, $$\phi_*B \cong \bigoplus_{u=1}^e \OO_{\PP^n} (-d_u), \ 0 \leq d_1 \leq \cdots \leq d_e$$ and $d_1 < d_e$. 

\begin{notation}\label{not-partialflag}
    Let $V = H^0(Z, B|_Z)$ and let $U= H^0(Z, \OO_Z)$. Both $U$ and $V$ are $e$-dimensional vector spaces. Multiplication by $x_0^{-v}$ identifies $B(vL)|_Z$ with $B|_Z$. For every integer $v \geq 0$ define
$$F_v = \mbox{im} \left(H^0(X, B(vL)) \to V \right).$$
Multiplication by $x_0$ gives an inclusion $F_v \subseteq F_{v+1}$, so that we get a partial flag. Evaluating at the point $[1:0:\cdots :0]$ gives
$$\dim (F_v) = \# \{ u \ | \ d_u \leq v \}.$$ In particular, $F_v = V$ for $v \geq d_e$ and since $d_1 < d_e$,  we have $0 < \dim F_{d_1} < e$. Hence, the partial flag $F_{\bullet}$ is non-trivial.

For a subset $T \subseteq Z$ of points of cardinality $|T|$, set 
$$V_T =H^0(Z, B|_Z \otimes \mathcal{I}_{Z \setminus T, Z})= H^0(T, B|_T) \quad \mbox{and} \quad F_v(T) = F_v \cap V_T.$$ Thus $F_v(T)$ consists of the elements of the evaluation image $F_v$ that vanish on $Z\setminus T$.
\end{notation}

Since sections in $F_v(T)$ are sections in $F_v$ that are required to vanish on $e-|T|$ points in the complement of $T$, we expect the dimension to be $\max(0, \dim(F_v) - e + |T|)$. Our first lemma shows that for a general choice of $Z$ this is indeed the case.

\begin{lemma}\label{lem-uniformposition}
  The section $Z$ can be chosen so that for every $v$ and every $T \subseteq Z$, $$\dim (F_v(T)) = \max(0, \dim(F_v) - e + |T|).$$
  Furthermore, this choice can be made simultaneously for all $B_i$ with $0 \leq i \leq n$.
\end{lemma}
\begin{proof}
   By Bertini's theorem, let $C$ be the integral projective curve defined by the intersection of $n-1$ general members of $|L|$. Consider the linear series $\tilde{U}_v$ on $C$ arising as the image of the restriction map
   $$\tilde{U}_v := \mbox{im} \left(H^0(X, B(vL)) \to H^0(C, B(vL)|_C) \right).$$
   Let $\PP^{n_0}$ be the linear span of $C$, and let $U\subset (\PP^{n_0})^*$ be the open set of hyperplanes meeting $C$ in $e$ distinct smooth points. For an integral complex projective curve, the monodromy over $U$ is the full symmetric group $\mathfrak S_e$ \cite[Proposition~1.1]{Kadets} (see also \cite[pp.~85--86]{harris} and \cite[\S 1]{PS}). For each $1\le k\le e$, this group acts transitively on the $k$-element subsets, so the finite \'etale cover of $U$ parametrizing such subsets of the hyperplane section is irreducible.

   On this cover the evaluation map has constant source $\tilde U_v$ and target the vector bundle with fiber $H^0(T',B(vL)|_{T'})$ at a $k$-element subset $T'$. Its rank drops below its generic value on a proper closed subset. The image of this subset in $U$ is closed and of smaller dimension because the cover is finite. Away from that image, every $k$-element subset of a given hyperplane section has the same evaluation rank. Only finitely many $v\ge 0$ need to be considered. Indeed, for $v\gg0$, Serre vanishing makes restriction from $X$ onto $H^0(C,B(vL)|_C)$ surjective and gives $H^1(C,B((v-1)L)|_C)=0$. Since $Z$ is a hyperplane divisor on $C$, evaluation onto $H^0(Z,B(vL)|_Z)$ is then surjective for every transverse $Z$. We can therefore choose a single general $Z$ for all $k$, all $v$, and all the finitely many $B_i$.

   After trivializing the fibers at the $e$ points, the coordinate functionals span $F_v^\vee$. A subset of $\dim F_v$ of them is a basis, so some $k$-element subset has evaluation rank $\min(k,\dim F_v)$. The common-rank assertion gives this rank for every such subset. Taking the kernel of evaluation on $Z\setminus T$ inside $F_v$ yields
   \[
     \dim F_v(T)=\dim F_v-\min(e-|T|,\dim F_v)=\max(0,\dim F_v-e+|T|),
   \]
   including the case $T=Z$, when the evaluation target is zero.
\end{proof}

\begin{notation}\label{not-r}
    Let $r$ be the least positive integer $v$ such that the evaluation maps $$H^0(X, vL) \to H^0(Z,vL|_Z) \quad \mbox{and} \quad H^0(X, B_i(vL)) \to H^0(Z, B_i(vL)|_Z)$$ are surjective for every $0 \leq i \leq n$. Multiplication by $x_0$ shows that these maps remain surjective for $v > r$.
\end{notation}

\begin{definition}
For every nonempty subset $T \subseteq Z$, define the \defi{weight} of $T$ by
$$w(T) = \sum_{v\geq 0} (|T| - \dim(F_v(T))).$$
Since $\dim(F_v(T))= |T|$ for $v \geq r$, this is a finite sum.
\end{definition}

\begin{lemma}\label{lem-weightcomparison}
    If $Z$ is chosen to satisfy Lemma \ref{lem-uniformposition}, then for every nonempty proper subset $T \subsetneq Z$, we have 
    $$w(T) = \sum_{k=1}^{|T|} d_{e-|T|+k} \quad \mbox{and} \quad \frac{w(T)}{|T|} > \frac{w(Z)}{e}.$$
\end{lemma}
 
\begin{proof}
We have  $\dim(F_v) = \#\{u|d_u \leq v\}$ and $\dim(F_v(T))= \max(0, \dim(F_v) - e + |T|)$. For each of the $|T|$ largest $d_u$, for $0 \leq v \leq d_u -1$, we get a contribution of $1$ to the sum defining $w(T)$. Hence, 
$$w(T) =  \sum_{v\geq 0} (|T| - \dim(F_v(T))) =  \sum_{k=1}^{|T|} d_{e-|T|+k}.$$ Hence $w(T)/|T|$ is the average of the $|T|$ largest $d_u$, whereas $w(Z)/e$ is the average of all the $d_u$ occurring in the splitting type of $\phi_*B$. Since $d_1 < d_e$, the first average is strictly larger.
\end{proof}

\begin{notation}
   For a positive integer $\ell$, let 
$$Q_{\ell} = \bigoplus_{v\geq 0} F_v/F_{v-\ell} \quad \mbox{and} \quad Q_{\ell}(T) = \bigoplus_{v\geq 0} F_v(T)/F_{v-\ell}(T),$$ where $F_v =0$ for $v<0$. Since $F_v = V$ and $F_v(T)= V_T$ for large enough $v$, both $Q_{\ell}$ and $Q_{\ell}(T)$ are finite-dimensional. Since $$(F_v \cap V_T) \cap F_{v-\ell} = F_{v-\ell} \cap V_T,$$ we get natural inclusions
$$\frac{F_v \cap V_T}{F_{v-\ell} \cap V_T} \hookrightarrow \frac{F_v}{F_{v-\ell}},$$ which induce the inclusion $Q_{\ell}(T) \hookrightarrow Q_{\ell}$.

Given a section $\sigma \in H^0(X, dL)$, the normalized restriction $\sigma/x_0^d|_Z$ multiplies $F_v$ into $F_{v+d}$ and $F_{v-\ell}$ into $F_{v-\ell+d}$. Hence, it induces multiplication maps $$m_{\sigma, v}: Q_{\ell, v} \to Q_{\ell, v+d}.$$
\end{notation}

\begin{definition}
    A \defi{subsystem preserved by multiplication} is a graded subspace $J = \bigoplus_v J_v \subseteq Q_{\ell}$ such that 
    $$m_{\sigma, v} (J_v) \subseteq J_{v+d}$$ for every $d \geq 0$, every $\sigma \in H^0(X, dL),$ and every $v$.
    For a nonzero finite-dimensional graded vector space $W$, define the \defi{mean degree} by 
    $$\mu(W) = \frac{\sum_v v \dim(W_v)}{\dim(W)}.$$
\end{definition}

\begin{lemma}\label{lem-computemuql}
    If $T$ is nonempty, then
    $$\dim(Q_{\ell}(T)) = \ell |T| \quad \mbox{and} \quad \mu(Q_{\ell}(T))= \frac{\ell-1}{2} + \frac{w(T)}{|T|}.$$
\end{lemma}

\begin{proof}
    Choose a basis of $H^0(T, B|_T)$ adapted to the filtration $F_{\bullet}(T)$. Assume the basis elements occur first in degrees $a_1, \dots, a_{|T|}$. As in the proof of Lemma \ref{lem-weightcomparison}, we have $w(T) = \sum_{k=1}^{|T|} a_k$.  In $Q_{\ell}(T)$, the $k$th basis element contributes one to the dimension in degrees $a_k, a_{k}+1, \dots, a_{k}+\ell-1$. Thus, the dimension of $Q_{\ell}(T)$ is $\ell |T|$ and the degree is
    $$\sum_{k=1}^{|T|} \sum_{h=0}^{\ell-1} (a_k + h) = \ell w(T) + \frac{|T|\ell (\ell-1)}{2}.$$ Dividing the degree by the dimension $\ell |T|$ gives the desired formula.
\end{proof}

We are now ready to prove the main result of this section.

\begin{proposition}\label{prop-multiplication}
    For all sufficiently large $\ell$, every proper nonzero  subsystem $J \subset Q_{\ell}$ preserved by multiplication satisfies 
    $$\mu(J) > \mu(Q_{\ell}).$$ Moreover, the integer $\ell$ can be chosen uniformly for finitely many evaluation filtrations on the same  $Z$.
\end{proposition}
\begin{proof}
    By Lemma \ref{lem-computemuql}, $$\mu(Q_{\ell}) = \frac{\ell-1}{2} + \frac{w(Z)}{e}.$$
    For a graded subspace $W$, set $$\Xi(W) := \sum_v \bigl(v-\mu(Q_{\ell})\bigr)\dim(W_v).$$
    For $W\ne0$, we have $\Xi(W)=\dim(W)\bigl(\mu(W)-\mu(Q_{\ell})\bigr)$. Hence, $\mu(J) > \mu(Q_{\ell})$ if and only if $\Xi(J) > 0$.
    
    Introducing $\Xi$ makes the calculation more transparent, because
     a vector in degree $v < \mu(Q_{\ell})$ contributes negatively to $\Xi(J)$ and a vector in degree $v > 
    \mu(Q_{\ell})$ contributes positively to $\Xi(J)$.

    By Notation~\ref{not-r}, evaluation on $Z$ is surjective for both $rL$ and $B(rL)$. Surjectivity for $B(rL)$ gives $F_r=V$.
    Let $c= \lfloor \mu(Q_{\ell}) \rfloor$.  For large $\ell$, the integer $c$ is approximately $\ell/2$. Consequently, for sufficiently large $\ell$, we have  $r < c$ and $r+c < \ell$.
Consequently,  $Q_{\ell, c} = V$. 

Let $T \subseteq Z$ be the support of $J_c$. First, assume that $T$ is a nonempty proper subset of $Z$. For $ v \leq c$, multiplication by $x_0^{c-v}$ is injective. Since the image of $J_v$ lies in $J_c$, its sections are supported on $T$. Hence, $$J_v \subseteq Q_{\ell}(T)_v \quad \mbox{for} \ v \leq c.$$ Since the weights $v - \mu(Q_{\ell})$ are nonpositive for these, the contribution of $J_v$ to $\Xi$ is at worst as negative as that of $Q_{\ell}(T)$ when $v \leq c$.

For $v\ge c+r$, the normalized restrictions of the sections in $x_0^{v-c-r}H^0(X,rL)$ give all functions on $Z$. Since $T$ is the support of $J_c$, we have $H^0(Z,\OO_Z)\cdot J_c=V_T$. Thus the products of these sections with $J_c$ span the image of $V_T$ in $Q_{\ell,v}=V/F_{v-\ell}$. This image is $Q_{\ell}(T)_v$, so $Q_{\ell}(T)_v\subseteq J_v$. In each degree $v\ge c+r$, the contribution of $J_v$ to $\Xi$ is therefore at least that of $Q_{\ell}(T)_v$.
    
For $c+1\le v\le c+r-1$, each graded piece has dimension at most $e$ and $0<v-\mu(Q_{\ell})\le v-c$. The total contribution lost in these degrees is at most $$E:=e\sum_{h=1}^{r-1}h=e\frac{r(r-1)}2,$$ independently of $\ell$. Thus $$\Xi(J) \geq \Xi(Q_{\ell}(T)) - E = \ell\left(w(T) - \frac{|T|}{e}w(Z) \right) - E.$$ By Lemma \ref{lem-weightcomparison}, the coefficient of $\ell$ is positive, so $\Xi(J)>0$ for $\ell\gg0$.

    If $T= \emptyset$, then $J_c=0$. It follows that $J_v =0$ for $v \leq c$. In this case, $\Xi(J)>0$, trivially.

    If $T=Z$, then $Q_{\ell}(T)= Q_{\ell}.$ By the choice of $r$ in Notation~\ref{not-r}, $F_v=V$ for $v\ge r$. Hence multiplication by $x_0$ maps $Q_{\ell,v}$ onto $Q_{\ell,v+1}$ for $v\ge r$, so $Q_{\ell}$ is generated in degrees at most $r$.  If $J_v = Q_{\ell,v}$ for $v \leq r$, then $J= Q_{\ell}$ contrary to the assumption that $J$ is proper. Hence, $J$ must miss a vector for $v \leq r < \ell/2$. Removing one such vector increases $\Xi$ by at least $\mu(Q_{\ell})-r$. The uncontrolled middle contributes at most $E$, so
    $$\Xi(J) \geq \mu(Q_{\ell})-r-E.$$ Since $\mu(Q_{\ell})$ goes to infinity with $\ell$, we conclude that $\Xi(J) > 0$ if $\ell \gg 0$. 
\end{proof}

\begin{remark}
Since $w(T) - \frac{|T|}{e} w(Z) \in \frac{1}{e} \mathbb{Z}_{>0}$, 
    a closer examination of the proof shows that choosing $$\ell > \max(4r, 2r+2E+1, eE)$$ suffices in Proposition \ref{prop-multiplication}.
\end{remark}

\subsection{The covering} Now we describe the covering $Y$ we will use throughout the proof.

\begin{notation}
  Let $(x_0, \dots, x_n)$ be general sections of $H^0(X, L)$ such that the common locus $V(x_1, \dots, x_n)= Z$ satisfies the conclusion of Lemma \ref{lem-uniformposition} for all $B_i$, $0 \leq i \leq n$. Let $\ell$ be an integer guaranteed by Proposition \ref{prop-multiplication} so that every proper nonzero subsystem $J$ of $Q_{\ell}$ preserved by multiplication satisfies $\mu(J) > \mu(Q_{\ell})$ for all $B_i$, $0 \leq i \leq n$.
\end{notation}

\begin{definition}
   Let $s= x_0^{\ell} + \sum_{i=1}^n \gamma_i x_i^{\ell}$ be a section of $H^0(X, L^{\otimes \ell})$, where $\gamma_i$ are general complex numbers.  Let $f:Y \to X$ be the cyclic cover of $X$ of degree $\ell$ branched along $V(s)$ defined by the equation $z^{\ell} = s$. Let $\pi:Y \to \PP^n$ be the projection given by $(x_1, \dots, x_n, z)$.
\end{definition}

\begin{remark}
 The sections $x_0,\ldots,x_n$ have no common zero. Where $x_1=\cdots=x_n=0$, we have $z^\ell=x_0^\ell\ne 0$, so $\pi$ is a morphism. Since $\pi^*\OO_{\PP^n}(1)=L_Y$ is ample, $\pi$ is finite. Locally $Y$ is cut out in $X\times\mathbb A^1$ by the monic equation $z^\ell-s=0$, so it is Cohen--Macaulay of pure dimension $n$. Thus $\pi$ is flat by miracle flatness.
 Moreover, modulo $x_1, \dots, x_n$, the cover $z^{\ell} = s$ becomes $z^{\ell} = x_0^{\ell}$ and geometrically realizes the truncations we have studied earlier in this section.  
\end{remark}

\begin{notation}
     Set $$S= \mathbb{C}[x_0, \dots, x_n]  \quad \mbox{and} \quad R= \bigoplus_v H^0(X, vL).$$
     Then $S$ is a polynomial subring of $R$ and $M_i \cong \bigoplus_{u=1}^e S(-d_{iu})$ is the corresponding module of twisted global sections of $\phi_* B_i$. Let $$N_i = M_i/(x_1, \dots, x_n)M_i.$$ Then $N_i \subset \mathbb{C}^e [x_0]$ with $(N_i)_v = x_0^v F_{iv}$. Let $$S' = \mathbb{C}[x_1, \dots, x_n , z], \quad R' = R[z]/(z^{\ell}-s), \quad M_i' = M_i\otimes_R R'$$
     $$\mathcal{T}= R'/S_+'R', \quad Q_i = M_i'/S_+' M_i' = N_i /x_0^{\ell} N_i.$$
\end{notation}

\begin{remark}
To compare the two polynomial rings, write
\[
 P=S[z]/(z^\ell-s)
   \cong S'[x_0]\big/\bigl(x_0^\ell-z^\ell+\textstyle\sum_{j=1}^n\gamma_j x_j^\ell\bigr).
\]
The equation on the right is monic in $x_0$, so $1,x_0,\ldots,x_0^{\ell-1}$ form an $S'$-basis of $P$. Since $R'=R\otimes_S P$, we have $M_i'=M_i\otimes_S P$. If $(b_u)_{u=1}^e$ is a homogeneous $S$-basis of $M_i$, then the elements $x_0^h b_u$, with $0\le h<\ell$ and $1\le u\le e$, form a homogeneous $S'$-basis of $M_i'$. This determines the splitting type
$$\pi_* A_i = \pi_* f^* B_i \cong \bigoplus_{h=0}^{\ell -1} \bigoplus_{u=1}^e \OO_{\PP^n}(-d_{iu} - h).$$
Modulo $x_1,\ldots,x_n,z$, the cover equation gives $x_0^\ell=0$. Since $(N_i)_v=x_0^v F_{iv}$, it follows that
$$ (Q_i)_v \cong (N_i)_v\big/x_0^\ell (N_i)_{v-\ell}
   \cong F_{iv}/F_{i,v-\ell}. $$
Under this identification, the class in $\mathcal T_d$ of a section $\sigma\in H^0(X,dL)$ acts from degree $v$ to degree $v+d$ by multiplication by $(\sigma/x_0^d)|_Z$. Thus each graded $\mathcal T$-submodule of $Q_i$ is a subsystem preserved by multiplication, and Proposition \ref{prop-multiplication} applies to every proper nonzero such submodule.
\end{remark}

We have therefore achieved our two goals on $Y$:
\begin{enumerate}
    \item The bundles $\pi_* A_i$ are explicit direct sums of line bundles, so that their cohomology is easy to compute.
    \item The multiplicity spaces of $\pi_* A_i$ form 
graded $\mathcal{T}$-modules in which every proper nonzero submodule satisfies the inequality of Proposition \ref{prop-multiplication}.
\end{enumerate}

\section{Constructing a bundle $G$ on $\PP^n$ so that the Hilbert polynomial has the desired form}\label{sec-numerical}
We preserve the notation from \S \ref{sec-prelim}. By Lemma \ref{lem-passtocover}, to construct an Ulrich bundle on $X$ for $mH$, it suffices to construct an Ulrich bundle on the cyclic cover $Y$ for $f^*\OO_X(mH)$. Set $D= f^*H$ and $L_Y= f^*(qH)$. The goal is to construct the Ulrich bundle on $Y$ as the kernel of a map of bundles
\begin{equation}\label{eq-defK}
   0 \to K \to F^{\oplus 2N} \to F(\delta L_Y)^{\oplus N} \to 0, 
\end{equation}
where $$F = A_0 \otimes \pi^*G$$ for a vector bundle $G$ on $\PP^n$.

\begin{notation}
    Let $P_G(t)$ be the Hilbert polynomial of $G$ and let $$p_G(t) = \frac{n!}{\rk(G)} P_G(t)$$ be the monic polynomial associated to $P_G(t)$.
\end{notation}

To get an Ulrich bundle for $mD$, we need $$H^j(Y, K(-imD))= 0 \quad \mbox{for} \ 1\leq i \leq n \ \mbox{and all} \ j.$$  In particular, 
$$\chi(Y, K(-imD)) = 0 \quad \mbox{for} \ 1 \leq i \leq n.$$

\subsection{Solving for the Hilbert polynomial $P_G$} We first solve for the invariants of $G$ needed to achieve these vanishings and then construct a vector bundle $G$ with the desired Hilbert polynomial as a direct sum of homogeneous bundles on $\PP^n$.

\begin{lemma}\label{lem-hilbaitensorw}
    The Hilbert polynomial of $\pi_* A_i \otimes_{\OO_{\PP^n}} G$ is 
    $\sum_u b_{iu} P_G(t-u)$.
\end{lemma}

\begin{proof}
    By Lemma \ref{lem-directsumoflb} and Notation \ref{not-directsum}, $\pi_* A_i = \bigoplus_u \OO_{\PP^n}(-u)^{b_{iu}}$. Consequently, $$\pi_* A_i \otimes G \cong \bigoplus_u G(-u)^{\oplus b_{iu}}.$$
    Since the Hilbert polynomial is additive over direct sums and the Hilbert polynomial of $G(-u)$ is $P_G(t-u)$, we conclude that the Hilbert polynomial of $\pi_* A_i \otimes G$ is $\sum_u b_{iu} P_G(t-u)$ as claimed.
\end{proof}

\begin{definition}\label{def-opai}
    For a polynomial $Q(t)$ of degree at most $n$, define the averaging operator $\avgop_i$ by
    $$(\avgop_i Q)(t) := \frac{1}{e \ell} \sum_{u} b_{iu} Q(t-u).$$
\end{definition}

\begin{remark}
  By Lemma \ref{lem-hilbaitensorw}, the Hilbert polynomial of $\pi_*A_i \otimes_{\OO_{\PP^n}} G$ is given by $e\ell (\avgop_i P_G)(t)$.
\end{remark}

By Taylor's Theorem, a polynomial of degree at most $n$ has the expansion
$$Q(t-u) = \sum_{k=0}^n \frac{(-u)^k}{k!} \frac{d^kQ}{dt^k} (t).$$
Setting $\mathcal{D}= \frac{d}{dt}$, we can formally express the latter expression as $$Q(t-u) = e^{-u\mathcal{D}} Q(t).$$ Here the exponential $e^{-u\mathcal{D}}$ is the finite sum $\sum_{k=0}^n \frac{(-u)^k\mathcal{D}^k}{k!}$ since $\mathcal{D}^k$ vanishes on polynomials of degree at most $n$ for $k \geq n+1$. Observe also that the operators $\mathcal{D}$ and $\avgop_i$ commute.

\begin{lemma}\label{lem-expAi}
    On polynomials of degree at most $n$, the operators $\avgop_i$ satisfy the relation
    $$\avgop_i = \avgop_0 e^{-ia\mathcal{D}/q}.$$
\end{lemma}

\begin{proof}
We can check the equality of two operators by checking equality on a basis of the space of polynomials of degree at most $n$.   Set $\Phi_i(t) := \chi(Y, A_i(tL_Y)).$ Since $$f_*A_i \cong \bigoplus_{h=0}^{\ell -1} \OO_X((t_a-ia-qh)H),$$ we have
$$\Phi_i(t)= \sum_{h=0}^{\ell -1} p(t_a - ia + q(t-h)),$$ where $p(t)=\chi(X,\OO_X(tH))$. Hence, we have the relation $$\Phi_i(t) = \Phi_0(t - ia/q).$$ Set $\zeta_n(t) = \binom{t+n}{n}$. The fact that $\pi_* A_i$ splits into a direct sum of line bundles implies that  $$\Phi_i (t)= e \ell \avgop_i \zeta_n(t).$$ Therefore, $\avgop_i$ and $\avgop_0 e^{-ia\mathcal{D}/q}$ agree on $\zeta_n(t)$. Since $\mathcal{D}$ and $\avgop_i$ commute, we conclude that they agree on all the derivatives of $\zeta_n(t)$, as well. The polynomials $\zeta_n, \mathcal{D}\zeta_n, \dots, \mathcal{D}^n \zeta_n$ are all nonzero and have strictly decreasing degree, hence they form a basis of the vector space of polynomials of degree at most $n$. Since the two operators $\avgop_i$ and $\avgop_0 e^{-ia\mathcal{D}/q}$ agree on this basis, they   must be equal.
\end{proof}

Comparing the coefficients of $\mathcal{D}$ in the expansion of $\avgop_i$ in Lemma \ref{lem-expAi}, we deduce the following corollary.

\begin{corollary}
    Set $$\mu_i = \frac{1}{e\ell} \sum_u u b_{iu}.$$ Then $$\mu_i = \mu_0 + \frac{ai}{q}.$$
\end{corollary}

\begin{notation}
  Set $u_{\max} = \max\{u \mid b_{iu}\ne 0 \text{ for some } 0\le i\le n\}$. Choose an integer $\delta\ge u_{\max}+5$, fixed as $M$ varies. By the division algorithm, write $m = M q + a$ with $0 \leq a < q$. Set
  $$\rho = \frac{m}{q} = M + \frac{a}{q}.$$
\end{notation}

If $K$ is defined by the exact sequence \eqref{eq-defK}, then the condition
$$\chi(K(-imD)) =0 \quad \mbox{for}  \ 1 \leq i \leq n$$  is equivalent to $$2 \chi(F(-imD)) - \chi(F(\delta L_Y-imD))=0 \quad  \mbox{for} \ 1 \leq i \leq n.$$ Since $F=A_0 \otimes_{\OO_Y} \pi^* G$, we have
 $$F(-imD) \cong A_0(-iaD) \otimes_{\OO_Y} (\pi^* G) (-iMqD) \cong A_i \otimes \pi^*(G (-iM)).$$ Hence, we would like to solve the equations 
$$2(\avgop_i P_G)(-iM)- (\avgop_i P_G)(\delta - iM)=0 \ \mbox{for} \ 1 \leq i \leq n.$$
By Lemma \ref{lem-expAi}, this is the same as solving the equation
$$2(\avgop_0 P_G) (-i\rho) - (\avgop_0 P_G)(\delta - i \rho)=0 \ \mbox{for} \ 1 \leq i \leq n.$$
In other words, the polynomial $$\avgop_0 (2 P_G(t) - P_G(\delta+t))$$ is a polynomial of degree $n$ with roots $-i \rho$ for $1 \leq i \leq n$, hence it must be a constant multiple of the following polynomial
$$\Pi_{\rho}(t) := \prod_{i=1}^n (t + i \rho).$$

\begin{notation}
    Set $\Delta_{\delta}$ to be the operator acting on polynomials of degree at most $n$ by 
$$\Delta_{\delta} Q(t) = Q(t+\delta) - Q(t).$$
\end{notation}

\begin{remark}
   The operator $\Delta_{\delta}$ cancels the leading coefficient of $Q(t)$, hence lowers the degree of the polynomial $Q(t)$. Therefore,  $\Delta_{\delta}^{n+1} = 0.$ Hence, the operator $1- \Delta_{\delta}$ is invertible with inverse
   $$(1- \Delta_{\delta})^{-1}= 1 + \Delta_{\delta} + \Delta_{\delta}^2 + \cdots + \Delta_{\delta}^n.$$
\end{remark}

\begin{remark}\label{rem-delta}
   The purpose of introducing $\Delta_{\delta}$ is the following useful relation
    $$(1-\Delta_{\delta})Q(t) = Q(t) - (Q(t+\delta) - Q(t))= 2 Q(t) - Q(t+\delta),$$ which occurs in the equation we are trying to solve.
\end{remark}

\begin{proposition} \label{prop-targethilb}
   If $K$ defined by the exact sequence \eqref{eq-defK} is an Ulrich bundle for $mD$ with $F= A_0 \otimes_{\OO_Y} \pi^* G$, then the   polynomial $p_G(t) = \frac{n!}{\rk(G)} P_G(t)$ is the polynomial
   $$p_G(t) = \avgop_0^{-1} (1 - \Delta_{\delta})^{-1} \Pi_{\rho}(t).$$
\end{proposition}

\begin{proof}
Observe that $\avgop_0$ and $\Delta_{\delta}$ commute. Moreover, since $1-\avgop_0$ lowers the degree of a polynomial it acts on, $(1-\avgop_0)^{n+1} = 0$, hence $\avgop_0$ is an invertible operator and we have a well-defined $\avgop_0^{-1}$.

   The polynomial  $p_G= \frac{n!}{\rk(G)} P_G(t)$ is a monic polynomial.  Observe that $\avgop_0$ and $1- \Delta_{\delta}$ yield a monic polynomial when applied to a monic polynomial. This is the reason to divide by the rank of $\pi_* A_0$ in the definition of the operator $\avgop_0$. By Remark \ref{rem-delta}, the equation we would like to solve can be expressed as
   $$\avgop_0 (1-\Delta_{\delta}) (p_G)(t) = \Pi_{\rho}(t)$$ up to a constant since both polynomials have the same roots. Since both polynomials are monic, the constant is one and equality holds.
   Hence, $$p_G (t) = \avgop_0^{-1} (1 - \Delta_{\delta})^{-1} \Pi_{\rho}(t)$$ as claimed, where we have used the fact that the operators commute.
\end{proof}

\begin{remark}
    For every $n$, $(1-\Delta_{\delta})^{-1} \Pi_{\rho}(t)$ is an explicit polynomial that can easily be computed. This is the main reason to swap the order of the operators using commutativity. For example when $n=2$, 
    \begin{eqnarray*}
        (1-\Delta_{\delta})^{-1} \Pi_{\rho}(t) &=& \Pi_{\rho} (t + 2\delta)- \Pi_{\rho}(t+ \delta) + \Pi_{\rho}(t) \\ &=& t^2 + (2\delta + 3 \rho)t + 3\delta^2 + 3 \delta \rho + 2 \rho^2
    \end{eqnarray*}
    
\end{remark}

Now that we determined the invariants of the bundle $G$ we need to construct, we construct the actual bundle. We will use direct sums of homogeneous bundles to construct $G$. We give a brief interlude recalling basic facts about homogeneous bundles on $\PP^n$.

\subsection{Homogeneous bundles on $\PP^n$} 
The Borel--Bott--Weil theorem computes the cohomology of homogeneous bundles on flag varieties \cite[Theorem~2.3]{CCHMW}. We use its projective-space form for Schur functors of the tautological quotient bundle \cite[Theorem~1.2, pp.~910--911]{EFW}.

Let $Q$ be the tautological quotient bundle on $\PP^n = \PP \Lambda$
$$0 \to \OO_{\PP^n}(-1) \to \OO_{\PP^n} \otimes \Lambda \to Q \to 0.$$ The Hilbert polynomial of $Q$ has roots at $t=-1, -2, \dots, -n+1, -n-1$.
\begin{definition}
    For a partition $\lambda_1 \geq \cdots \geq \lambda_n \geq 0$, let $G_{\lambda}$ be the homogeneous bundle $\mathbb{S}_{\lambda}(Q)$. 
\end{definition}

\begin{remark}
    The Hilbert polynomial of $G_{\lambda}$ has roots 
$$\beta_i = -i - \lambda_{n+1-i}  \ \mbox{for} \ 1 \leq i \leq n.$$
Hence,
$$\chi(G_{\lambda}(t)) = \frac{\rk(G_{\lambda})}{n!} \prod_{i=1}^n (t - \beta_i).$$ The rank is determined by the Weyl dimension formula, but will not play a role in our discussion.
\end{remark}

Set $\beta_0 = +\infty$ and $\beta_{n+1} = -\infty$. 

\begin{theorem}[Borel-Bott-Weil Theorem]\label{thm-BBW}
   If $t= \beta_i$ for some $i$, then $H^j(\PP^n, G_{\lambda}(t))=0$ for all $j$. If $\beta_h > t > \beta_{h+1}$ is an integer, then $H^j(\PP^n, G_{\lambda}(t))=0$ for all $j \not= h$ and $H^h(\PP^n, G_{\lambda}(t))$ is an irreducible $\SL(\Lambda)$-representation. 
\end{theorem}

Every strictly decreasing sequence of integers $\beta_1>\cdots>\beta_n$ with $\beta_1\le-1$ is realized by the partition $\lambda_{n+1-i}=-i-\beta_i$. Twisting a bundle by $\OO_{\PP^n}(c)$ shifts all its roots by $-c$. Thus every strictly decreasing sequence of $n$ integers is realized by a line-bundle twist of a Schur functor of $Q$, as in \cite[Theorem~6.2, p.~879]{ES09}.

\subsection{The construction of the bundle with the required Hilbert polynomial} We are now ready to construct $G$ as a direct sum of homogeneous bundles on $\PP^n$.

\begin{theorem}\label{thm-constructbundle}
    For all sufficiently large $M$, there exists a direct sum 
    $$G = \bigoplus_{\epsilon} G_{\epsilon}^{\oplus \nu_{\epsilon}},$$ where $G_{\epsilon}$ is a homogeneous bundle on $\PP^n$, $\nu_{\epsilon}$ is a positive integer, $$\chi(G(t))= \frac{\rk(G)}{n!} p_G(t),$$ where $p_G$ is the polynomial determined in Proposition \ref{prop-targethilb} and the roots of $G_{\epsilon}$ are all of the form 
    $$\beta_{i, \epsilon} = - i \rho - \delta - \mu_0 + e_{i, \epsilon}, \ \ \mbox{with} \ |e_{i, \epsilon}| \leq 3.$$
    Moreover, for any two indices $\epsilon, \epsilon'$, $|\beta_{i, \epsilon} - \beta_{i, \epsilon'}|\leq 2$.
\end{theorem}

\begin{proof} 
Recall that $\rho = M + \frac{a}{q}.$

    {\bf Step 1.} We first show that if $M \gg 0$, then the roots of the target polynomial $p_G$ are real and far apart. By Proposition \ref{prop-targethilb}, $$p_G = \avgop_0^{-1}(1-\Delta_{\delta})^{-1} \Pi_{\rho}(t).$$ The operators $\avgop_0^{-1}$ and $(1- \Delta_{\delta})^{-1}$ are independent of $M$ and depend only on $\delta$ and $A_0$. Express
    $$ \avgop_0^{-1}(1-\Delta_{\delta})^{-1}  = 1 + (\mu_0 + \delta) \mathcal{D} + c_2 \mathcal{D}^2 + \cdots + c_n \mathcal{D}^n.$$ The higher order terms will not play a role in the argument.
    
    Make a change of variables $t = \rho u$. Then 
    $$\Pi_{\rho}(t) = \Pi_{\rho}(\rho u) = \prod_{i=1}^n (\rho u + i \rho) = \rho^n \prod_{i=1}^n (u+i).$$
    Set $\Psi(u) = \prod_{i=1}^n (u+i).$
    Hence, $$\rho^{-n} p_G(\rho u) = \Psi(u) + (\delta + \mu_0) \frac{1}{\rho} \Psi'(u) + \mbox{h.o.t.}$$
    When $\frac{1}{\rho} =0$, we find that the polynomial has simple roots at $-1, \dots, -n$. Hence, by the implicit function theorem, if $\rho \gg 0$, then the polynomial has real roots 
    $$-i - (\mu_0 + \delta) \frac{1}{\rho} + O((1/\rho)^2).$$
    If $\rho \gg 0$, these are $n$ distinct real numbers. Since the degree of the polynomial is $n$, we conclude that the roots are all real and distinct. Reversing the change of variables, we conclude that $p_G$ has distinct real roots of the form
    $$-i\rho - \mu_0 - \delta + O(1/\rho).$$ Hence, the roots are all real and approximately separated by $\rho$.

    {\bf Step 2.} Next, we express $p_G$ as a real positive linear combination of reduced Hilbert polynomials of homogeneous bundles.  When $\rho$ is sufficiently large, for each root $\alpha_i$ of $p_G$, choose two integers $r_{i,0} < \alpha_i < r_{i,1}$ so that $r_{i,1} - r_{i,0} \leq 2$. If $\alpha_i$ is an integer, let $r_{i,0}$ and $r_{i,1}$ be its two neighbors. Otherwise, let $r_{i,0}$ be the floor of $\alpha_i$ and let $r_{i,1}$ be the ceiling, so that $r_{i,1} - r_{i,0} =1$. Since the distance between the roots is approximately $\rho$, for large $M$ these distances are large, and these intervals do not overlap. We can thus assume that 
    $$r_{n,0} < r_{n,1} <  \cdots < r_{1,0} < r_{1,1}.$$ Since $\alpha_i \in (r_{i,0}, r_{i,1})$,
we can write $$(t-\alpha_i) = \theta_i (t -r_{i,0}) + (1-\theta_i) (t-r_{i,1}), \quad 0 < \theta_i <1.$$
We get $p_G = \prod_{i=1}^n (\theta_i (t -r_{i,0}) + (1-\theta_i) (t-r_{i,1})).$
For $M$ sufficiently large, we have $r_{1,1}\le-1$. For each $\epsilon=(\epsilon_1,\ldots,\epsilon_n)\in\{0,1\}^n$, take the partition $\lambda(\epsilon)_{n+1-i}=-i-r_{i,\epsilon_i}$ and set $G_\epsilon=\mathbb S_{\lambda(\epsilon)}Q$. Its roots are $r_{i,\epsilon_i}$ for $1\le i\le n$. Expanding the product expresses $p_G$ as a positive real linear combination of the reduced Hilbert polynomials of these $2^n$ homogeneous bundles.

{\bf Step 3.} Next, we replace the positive real coefficients with positive rational coefficients.  Consider the equation $$\sum_{\epsilon} c_{\epsilon} p_{G_{\epsilon}}(t) = p_G(t),$$ where we treat the coefficients $c_{\epsilon}$ as variables. Since the Hilbert polynomials are rational polynomials, this is a system with rational coefficients. By Gaussian elimination, since it has a real solution, it must also have a rational solution. Furthermore, the rational solutions are dense in the space of solutions. Let $w_{\epsilon}$ be a solution close to the real solution $c_{\epsilon}$ constructed in Step 2. Since $c_{\epsilon}$ are positive real numbers, we can arrange so that $w_{\epsilon}$ are positive rational numbers. Moreover, we must have $\sum_{\epsilon} w_{\epsilon} =1$ since $p_G$ is monic. Let $\nu_{\epsilon} = N \frac{w_{\epsilon}}{\rk(G_{\epsilon})}$ for an integer $N$ that clears the denominators of all the ratios $w_{\epsilon}/\rk(G_{\epsilon})$.

Setting $G = \bigoplus_{\epsilon} G_{\epsilon}^{\oplus \nu_{\epsilon}}$, we get a bundle with the desired  Hilbert polynomial
$\frac{\rk(G)}{n!} p_G(t)$. 

Finally, the $i$th root of $G_{\epsilon}$ is within an interval of $2$ of $\alpha_i$. If $M$ is sufficiently large, $\alpha_i$ is close to $-i\rho - \delta - \mu_0$, so the roots are within at most 3 of this number. This concludes the proof of the theorem.
    
\end{proof}

\section{Constructing surjective bundle maps and isomorphisms on cohomology}\label{sec-cohomology}
We preserve the notation introduced in the previous sections.
In \S \ref{sec-numerical}, we constructed a bundle $G$ on $\PP^n$ such that if $F= A_0 \otimes_{\OO_Y} \pi^* G$ on $Y$ and there exists an exact sequence 
$$0 \to K \to F^{\oplus 2N} \to F(\delta L_Y)^{\oplus N} \to 0,$$ then $K$ has vanishing Euler characteristic at the twists $-im D$. To show that $K$ is an Ulrich bundle, we  need to verify that there is a surjective morphism  $F^{\oplus 2N} \to F(\delta L_Y)^{\oplus N}$ which induces isomorphisms on cohomology for the twists $-im D$ with $1 \leq i \leq n$.

\subsection{Surjectivity of the maps between homogeneous bundles and their cohomologies} Let $\bfe_1, \dots, \bfe_n$ denote the standard unit vectors of length $n$, and put $\lambda_0=+\infty$. Pieri's rule \cite[Theorem~1.1, p.~910]{EFW} gives
$$Q \otimes \mathbb{S}_{\lambda} Q = \bigoplus_{ 1\leq j \leq n, \\ \lambda_{j-1}>\lambda_j} \mathbb{S}_{\lambda + \bfe_j} Q,$$ where $\lambda + \bfe_j$ adds one to the $j$th part provided the resulting sequence is still a partition.  The operation of adding $\bfe_j$ decreases the $(n+1-j)$th root of the Hilbert polynomial of the homogeneous bundle  by one. 

We can represent partitions by their corresponding Young diagrams. We say $\lambda' \supset \lambda$ if the Young diagram of $\lambda'$ contains the Young diagram of $\lambda$. This  gives a partial ordering on partitions.

\begin{lemma}\label{lem-repeatedPieri}
  Let $E_{\lambda} = \mathbb{S}_{\lambda} Q(c)$ and let $E_{\lambda'} = \mathbb{S}_{\lambda'} Q(c)$ such that $\lambda' \supset \lambda$. Fix a twist $w$ and a cohomological degree $h$. Suppose there is a sequence of partitions from $\lambda$ to $\lambda'$ that adds one box at a time and every bundle $E_{\lambda''}$ during this process satisfies $$H^{h+1}(E_{\lambda''}(w-1))=0.$$
  Then the evaluation map $$\Hom(E_{\lambda}, E_{\lambda'}) \otimes H^h(E_{\lambda}(w)) \to H^h(E_{\lambda'}(w))$$ is surjective.
\end{lemma}

\begin{proof}
Let $\lambda''$ be a partition that appears in the sequence of partitions. Tensor the Euler sequence by $E_{\lambda''}$, twist by $w$ and apply Pieri's rule to obtain
    $$0 \to E_{\lambda''}(w-1) \to E_{\lambda''}(w) \otimes \Lambda \to \bigoplus_{1\leq j \leq n, \\ \lambda_{j-1}''>\lambda_j''} E_{\lambda'' + \bfe_j}(w) \to 0.$$ Since $H^{h+1}(E_{\lambda''}(w-1))=0,$ the long exact sequence associated to this sequence gives the surjectivity of 
    $$\Lambda \otimes H^h(E_{\lambda ''} (w)) \to \bigoplus_{1\leq j \leq n, \\ \lambda_{j-1}''>\lambda_j''} H^h(E_{\lambda'' + \bfe_j}(w)).$$ Hence, the projection to any $\lambda'''$ of the form $\lambda'' + \bfe_j$ remains surjective. By induction, we conclude that $\Lambda^{\otimes (|\lambda'|-|\lambda|)} \otimes H^h(E_{\lambda}(w))$ surjects onto $H^h(E_{\lambda'}(w))$. Since the pure tensors in $\Lambda^{\otimes (|\lambda'|-|\lambda|)}$ give homomorphisms from $E_{\lambda}$ to $E_{\lambda'}$  by composing the homomorphisms between the bundles along the path,  $\Lambda^{\otimes (|\lambda'|-|\lambda|)}$ maps to $\Hom(E_\lambda, E_{\lambda'})$. The lemma follows.
\end{proof}

We apply the lemma to the bundle $G$ we constructed in Theorem \ref{thm-constructbundle}.

\begin{proposition}\label{prop-surjhom}
   For all $1 \leq i \leq n$, all sufficiently large $M$ and all indices $\epsilon$ and $\eta$, all $s$ with $\mathcal{T}_s \not= 0$ and all $j$ for which $(Q_i)_j$ and $(Q_i)_{j+s}$ are nonzero, the evaluation map 
   $$\Hom(G_{\epsilon}, G_{\eta}(\delta -s))\otimes H^{i-1}(G_{\epsilon}(-iM-j)) \to H^{i-1}(G_{\eta}(-iM+\delta -j-s))$$
   is surjective. Furthermore, all the other cohomology groups of these bundles are zero.
\end{proposition}

\begin{proof}
Let $\mathbf 1=(1,\ldots,1)$ and put $\lambda'=\lambda(\eta)+(\delta-s)\mathbf 1$. The Euler sequence gives $\det Q\simeq(\det\Lambda)\otimes\OO_{\PP^n}(1)$. Trivializing the constant line $\det\Lambda$, we obtain
$$G_\eta(\delta-s)\cong\mathbb S_{\lambda'}Q.$$
By Theorem \ref{thm-constructbundle}, the corresponding roots of $G_\epsilon$ and $G_\eta$ differ by at most two. Since $j\ge0$ and $j+s\le u_{\max}$, we have
$$\lambda'_k-\lambda(\epsilon)_k
 =\delta-s+\beta_{n+1-k,\epsilon}-\beta_{n+1-k,\eta}
 \ge\delta-u_{\max}-2\ge3.$$
Thus $\lambda'$ contains $\lambda(\epsilon)$. Take the path from $\lambda(\epsilon)$ to $\lambda'$ obtained by completing each row in turn, starting with the first.

Set $w=-iM-j$. Each added box decreases one root by one. The $i$th root is therefore largest at the source, where
$$w-\beta_{i,\epsilon}=\delta+\mu_i-j-e_{i,\epsilon}\ge2$$
by $j\le u_{\max}$, $\delta\ge u_{\max}+5$, $\mu_i\ge0$ and $|e_{i,\epsilon}|\le3$. The $(i-1)$st root is smallest at the target. For $i>1$ its value there is $\beta_{i-1,\eta}-(\delta-s)$, which exceeds $w$ for large $M$ by Theorem \ref{thm-constructbundle}. For $i=1$ the upper endpoint is $+\infty$. Thus both $w$ and $w-1$ lie between the $i$th and $(i-1)$st roots at every step. The Borel-Bott-Weil Theorem gives concentration in degree $i-1$ at both twists, and Lemma \ref{lem-repeatedPieri} gives the asserted surjectivity.
\end{proof}

\subsection{Getting isomorphisms on the cohomology}
Write $E_\epsilon=\mathbb C^{\nu_\epsilon}$ and $F=A_0\otimes\pi^*G$. For fixed $i$ put
\[
 H_{\epsilon,j}=H^{i-1}(G_\epsilon(-iM-j)),\qquad
 H'_{\epsilon,j}=H^{i-1}(G_\epsilon(-iM+\delta-j)).
\]
Choose a homogeneous $S'$-basis of $M_i'$ as in \S\ref{sec-cover}. Its degree-$j$ elements give a basis of $(Q_i)_j$ and identify $\pi_*A_i$ with $\bigoplus_j\OO_{\PP^n}(-j)\otimes (Q_i)_j$. Since $\pi$ is finite, the projection formula gives
\[
\begin{aligned}
 U_i&:=H^{i-1}(Y,F(-imD))
      \cong\bigoplus_{\epsilon,j} H_{\epsilon,j}\otimes (Q_i)_j\otimes E_\epsilon,\\
 V_i&:=H^{i-1}(Y,F(\delta L_Y-imD))
      \cong\bigoplus_{\epsilon,j} H'_{\epsilon,j}\otimes (Q_i)_j\otimes E_\epsilon.
\end{aligned}
\]

\begin{lemma}
    $2\dim (U_i) = \dim (V_i)$
\end{lemma}

\begin{proof}
Since all the cohomology except the $(i-1)$st cohomology of $A_0 \otimes_{\OO_Y} \pi^*G (-imD)$ vanishes by construction, these dimensions are $(-1)^{i-1}$ times the Euler characteristic. Since $$2 \chi(A_0 \otimes \pi^*G(-imD))= \chi(A_0 \otimes \pi^*G(\delta L_Y -imD))$$ by construction, the lemma follows.
\end{proof}

We would like to construct an isomorphism $U_i^{\oplus 2N} \to V_i^{\oplus N}$ for sufficiently large $N$ using maps between homogeneous bundles discussed in the previous subsection, multiplication maps on $(Q_i)_j$ and linear maps between $\mathbb{C}^{\nu_\epsilon}$. On $\PP^n$, Coskun, Huizenga, and Smith showed that for fixed $0\le d<e$ and $s+t$ sufficiently large, a general map $\OO(d)^{\oplus s}\to\OO(e)^{\oplus t}$ induces a map of maximal rank on global sections \cite[Corollary~2.8]{CHS}. For the cohomology spaces $U_i$ and $V_i$, we use the following criterion for matrix spaces.

\begin{definition}
Let $U$ and $V$ be two vector spaces of the same dimension.   Let $\mathcal{B} \subset \Hom(U,V)$ be a subspace of homomorphisms and let $W \subset U$ be a subspace. Let $\mathcal{B} W$ be the subspace spanned by all elements of the form $b(w)$ for $b \in \mathcal{B}$ and $w \in W$. The subspace $W$ is \defi{shrunk}
 if $\dim (\mathcal{B} W) < \dim (W).$
 \end{definition}

The following criterion follows from the characterization of rank over the free skew field by Fortin and Reutenauer \cite[Theorem~1, p.~4]{FR} and the matrix evaluation criterion in \cite[Proposition~1.12, p.~48]{DM}.
\begin{theorem}\label{thm-isomscaling}
    A  matrix space $\mathcal{B} \otimes \operatorname{Mat}_{N \times N} \subset \Hom(U^{\oplus N}, V^{\oplus N})$ contains an invertible element for some $N$ if and only if $\mathcal{B}$ has no shrunk subspace.
\end{theorem}

\begin{corollary}\label{cor-noshrunk}
    Suppose $\dim(V)=2\dim(U)$ and $\mathcal B\subset\Hom(U,V)$ satisfies $\dim(\mathcal B W)\geq 2\dim(W)$ for every $W\subset U$. Then for every sufficiently large $N$, the space $\mathcal B\otimes\operatorname{Mat}_{N\times 2N}$ contains an isomorphism $U^{\oplus 2N}\to V^{\oplus N}$.
\end{corollary}

\begin{proof}
The assertion is immediate if $U=0$. Otherwise put $d=\dim V=2\dim U$ and set $\mathcal C=\mathcal B\otimes\Hom(\mathbb C^2,\mathbb C)\subset\Hom(U^{\oplus 2},V)$. For $W'\subset U^{\oplus 2}$, let $\widetilde W=p_1W'+p_2W'$, where $p_1,p_2$ are the coordinate projections. Then $\mathcal C W'=\mathcal B\widetilde W$ and $\dim W'\leq 2\dim\widetilde W$. Thus $\mathcal C$ has no shrunk subspace. By Theorem~\ref{thm-isomscaling}, $\mathcal C\otimes\operatorname{Mat}_{N\times N}$ contains an invertible matrix for some $N$. A tuple of $d\times d$ matrices spanning $\mathcal C$ therefore lies outside the null cone by \cite[Corollary~1.6, p.~47]{DM}. Since $d\ge2$, \cite[Theorem~1.8, p.~47]{DM} gives an invertible matrix in $\mathcal C\otimes\operatorname{Mat}_{N\times N}$ for every $N\ge d-1$. This is the required rectangular matrix space after identifying $(U^{\oplus 2})^{\oplus N}$ with $U^{\oplus 2N}$.
\end{proof}

For all the twists we use the same parameter space
\[
 \mathcal P=\bigoplus_{\epsilon,\eta,s}
 \Hom(G_\epsilon,G_\eta(\delta-s))\otimes\mathcal T_s\otimes\Hom(E_\epsilon,E_\eta).
\]
For $\alpha\otimes b\otimes C$ in the $(\epsilon,\eta,s)$ summand, define $\Gamma_i(\alpha\otimes b\otimes C):U_i\to V_i$ by its action on each $(\epsilon,j)$ summand,
\[
 x\otimes v\otimes u\longmapsto
 H^{i-1}(\alpha(-iM-j))(x)\otimes bv\otimes C(u)
 \quad\text{in }H'_{\eta,j+s}\otimes (Q_i)_{j+s}\otimes E_\eta.
\]
It is zero on source summands with index different from $\epsilon$. In particular the same $\alpha$, $b$ and $C$ act at every degree $j$. Let $\mathcal B_i$ be the image of the resulting linear map $\Gamma_i:\mathcal P\to\Hom(U_i,V_i)$.

\begin{proposition}\label{prop-noshrunk}
    For sufficiently large $M$, simultaneously for $1 \leq i \leq n$, 
    $$\dim \mathcal{B}_i W \geq 2 \dim(W), \quad \mbox{for every subspace} \ W \subset U_i.$$
\end{proposition}

\begin{proof}
    Define the defect by 
    $$ \defe(W) = 2 \dim(W) - \dim(\mathcal{B}_i W). $$
    The proposition is equivalent to proving that $\defe(W) \leq 0$ for every subspace $W \subset U_i$.

{\bf Step 1.} We first show that there is a largest subspace of maximal defect. For any two subspaces $W_1,W_2\subseteq U_i$, the spans of their images satisfy
\[
 \mathcal B_i(W_1+W_2)=\mathcal B_iW_1+\mathcal B_iW_2,\qquad
 \mathcal B_i(W_1\cap W_2)\subseteq\mathcal B_iW_1\cap\mathcal B_iW_2.
\]
Applying the dimension formula for a sum of subspaces gives

\begin{equation}\label{eq-defect}
    \defe(W_1) + \defe(W_2) \leq \defe(W_1 + W_2) + \defe(W_1 \cap W_2).
\end{equation}

Since the defect takes values that are finitely many integers, some subspace achieves the maximum value. If two subspaces achieve the maximum, by the Inequality \eqref{eq-defect}, both the sum and the intersection must also have the same defect. Hence, there is a largest subspace with the maximal defect. Write $W_{\max}$ for this subspace.

{\bf Step 2.} The group $\SL(\Lambda)$ acts on the homogeneous bundles and $(g_\epsilon)_\epsilon\in\prod_\epsilon\operatorname{GL}(E_\epsilon)$ acts on their multiplicity spaces. Conjugation sends the operator associated to $\alpha\otimes b\otimes C$ to one of the same form, by equivariance on cohomology and the change $C\mapsto g_\eta Cg_\epsilon^{-1}$. The grading torus acts by $t^j$ on $(Q_i)_j$, so conjugation multiplies this operator by $t^s$. Thus these groups preserve $\mathcal B_i$ and the defect, and hence preserve $W_{\max}$.

The grading characters distinguish the degrees $j$, and the scalar centers of the separate $\operatorname{GL}(E_\epsilon)$ distinguish the indices $\epsilon$. By Borel--Bott--Weil, $H_{\epsilon,j}$ is irreducible for $\SL(\Lambda)$ whenever $(Q_i)_j\ne0$. Hence $H_{\epsilon,j}\otimes E_\epsilon$ is irreducible for the product group. Complete reducibility and Schur's lemma give subspaces $J_{\epsilon,j}\subset (Q_i)_j$ such that
\[
 W_{\max}=\bigoplus_{\epsilon,j} H_{\epsilon,j}\otimes J_{\epsilon,j}\otimes E_\epsilon .
\]

{\bf Step 3.} Let $J\subset Q_i$ be the graded $\mathcal T$-submodule generated by all $J_{\epsilon,j}$, so that
\[
 J_k=\sum_{\epsilon,\,j+s=k}\mathcal T_s J_{\epsilon,j}.
\]
Define $U_J$ and $V_J$ by replacing $(Q_i)_j$ with $J_j$ in the decompositions of $U_i$ and $V_i$. Then $W_{\max}\subset U_J$ and $\mathcal B_iW_{\max}\subset V_J$.

For the reverse inclusion, fix $\epsilon,\eta,j,s$ and $b\in\mathcal T_s$ with $bJ_{\epsilon,j}\ne0$. Then both $(Q_i)_j$ and $(Q_i)_{j+s}$ are nonzero. Apply the operators $\Gamma_i(\alpha\otimes b\otimes C)$ to vectors of $W_{\max}$ supported in the single summand $(\epsilon,j)$. As $\alpha$ varies, the images on cohomology span $H'_{\eta,j+s}$ by Proposition~\ref{prop-surjhom}. As $C$ varies, the images of $E_\epsilon$ span $E_\eta$. Thus these images span
\[
 H'_{\eta,j+s}\otimes bJ_{\epsilon,j}\otimes E_\eta.
\]
At fixed $\eta,k$, summing over $\epsilon$, $j+s=k$ and $b\in\mathcal T_s$ gives $H'_{\eta,k}\otimes J_k\otimes E_\eta$. Consequently $\mathcal B_iW_{\max}=V_J$, and
\[
 \defe(W_{\max})\leq 2\dim(U_J)-\dim(V_J).
\]
For $J=0$ or $J=Q_i$ the right hand side is zero. We may therefore assume that $J$ is a proper nonzero graded submodule.

{\bf Step 4.} Using Proposition \ref{prop-multiplication}, we show that $\defe(W_{\max}) \leq 0$. Define the operator acting on polynomials of degree at most $n$ by
$$\avgop_J = \frac{1}{\dim J} \sum_j (\dim (J_j))e^{-j \mathcal{D}},$$ where recall that $\mathcal{D} = \frac{d}{dt}$.
Purity of cohomology and the relation $$\avgop_0 (2 p_G(t) - p_G(t+\delta)) = \Pi_{\rho}(t)$$ imply that
$$\frac{2 \dim(U_J) - \dim(V_J)}{\rk(G)} = \frac{\dim(J)}{n!} (-1)^{i-1} (\avgop_J \avgop_0^{-1} \Pi_{\rho})(-iM) =\frac{\dim(J)}{n!} (-1)^{i-1} (\avgop_J \avgop_i^{-1} \Pi_{\rho})(-i\rho).$$
Since $$\avgop_J \avgop_i^{-1} = 1 + (\mu_i - \mu(J)) \mathcal{D} + \mbox{h.o.t.}$$ and $\Pi_{\rho}(-i\rho)=0$, we conclude that
$$\frac{2 \dim(U_J) - \dim(V_J)}{\rk(G)} = \frac{\dim(J)}{n!} (\mu_i - \mu(J)) (i-1)! (n-i)! \rho^{n-1} + O(\rho^{n-2}).$$
By Proposition \ref{prop-multiplication}, $\mu(J)>\mu(Q_i)=\mu_i$, so the leading coefficient is strictly negative. For fixed $i$, the module $Q_i$ is independent of $M$, and $\avgop_J$ depends only on the dimension vector of $J$. There are finitely many such vectors and finitely many twists $i$. Thus one lower bound on $M$ makes the expression negative for every proper nonzero $J$ at every required twist.
\end{proof}

\begin{corollary}\label{cor-isom}
For $M$ sufficiently large, there is an integer $N_0$ such that for every $N\geq N_0$ and each $1\leq i\leq n$, the general element of $\mathcal B_i\otimes\operatorname{Mat}_{N\times2N}$ is an isomorphism $U_i^{\oplus 2N}\to V_i^{\oplus N}$.
\end{corollary}

\begin{proof}
Apply Corollary~\ref{cor-noshrunk} to each of the finitely many spaces in Proposition~\ref{prop-noshrunk} and take the maximum of the resulting bounds. Invertibility is open in each matrix space.
\end{proof}

\subsection{Lifting the cohomology isomorphisms to a surjective bundle map}

\begin{proposition}\label{prop-actualbundle}
    Set $F=A_0\otimes_{\OO_Y} \pi^* G$ constructed in \S \ref{sec-numerical}. For sufficiently large $N$, there exists a surjective bundle map 
    $$\varphi: F^{\oplus 2N} \to F(\delta L_Y)^{\oplus N}$$ such that $\varphi$ induces isomorphisms  
    $$H^{i-1}(Y, F(-imD))^{\oplus 2N} \to H^{i-1}(Y, F(\delta L_Y-imD))^{\oplus N}$$ for $1 \leq i \leq n$. In particular, the kernel of $\varphi$ is an Ulrich bundle on $Y$ with respect to $mD$.
\end{proposition}

\begin{proof}
Choose a graded linear section of $R'\twoheadrightarrow\mathcal T$ taking $1$ to $1$, and denote the lift of $b\in\mathcal T_s$ by $\widetilde b$. Multiplying $\pi^*\alpha$ by the section $\widetilde b$ of $sL_Y$ and applying $C$ on multiplicity spaces gives a linear map
\[
 \mathcal P\longrightarrow\Hom(F,F(\delta L_Y)).
\]
For $\mathcal P_N=\mathcal P\otimes\operatorname{Mat}_{N\times 2N}$ this gives bundle maps $\varphi_p:F^{\oplus2N}\to F(\delta L_Y)^{\oplus N}$.

Fix $i$. With respect to the chosen $S'$-basis of $M'_i$, the entry of multiplication by $\widetilde b$ from a generator of degree $j$ to one of degree $k$ is a homogeneous polynomial of degree $s+j-k$. Entries of negative degree are zero. Multiply each nonzero entry by $\tau^{s+j-k}$, combine it with the corresponding $\alpha$ and $C$, and take cohomology. This gives a polynomial family $\Upsilon_i(\tau,p):U_i^{\oplus 2N}\to V_i^{\oplus N}$, linear in $p\in\mathcal P_N$. At $\tau=1$ it is the map induced by $\varphi_p$ at the $i$th twist. At $\tau=0$ only the constant entries survive. Reduction modulo $S'_+$ gives exactly these entries, so the surviving matrix is multiplication by $b$ on $Q_i=M'_i/S'_+M'_i$. Thus
\[
 \Upsilon_i(0,-)=\Gamma_i\otimes\operatorname{id}_{\operatorname{Mat}_{N\times 2N}},
 \qquad
 \Upsilon_i(0,\mathcal P_N)=\mathcal B_i\otimes\operatorname{Mat}_{N\times 2N}.
\]

For $\tau\ne0$, let $\sigma_\tau$ act on the degree-$s$ parameter by $\tau^s$, and let $D_U(\tau),D_V(\tau)$ act by $\tau^j$ on each copy of the degree-$j$ source and target summands. The exponent $s+j-k$ gives
\[
 \Upsilon_i(\tau,p)
  =D_V(\tau)^{-1}\Upsilon_i(1,\sigma_\tau p)D_U(\tau).
\]
Choose $N\geq\max(N_0,n)$ with $N_0$ as in Corollary~\ref{cor-isom}. For each $i$ the equality of images at $\tau=0$ gives $p_i\in\mathcal P_N$ with $\det\Upsilon_i(0,p_i)\ne0$. This determinant remains nonzero for some $\tau\ne0$, so $\det\Upsilon_i(1,\sigma_\tau p_i)\ne0$. Hence the parameters inducing an isomorphism at the $i$th twist form a nonempty Zariski open subset of $\mathcal P_N$.

Surjectivity on $Y$ is also open in $\mathcal P_N$, since the rank-drop locus in $Y\times\mathcal P_N$ is closed and $Y$ is projective. To see that this open set is nonempty, take $b=1$ and use scalar maps on $G$. This yields the pullbacks of all $N\times2N$ matrices of sections of $\OO_{\PP^n}(\delta)$, tensored with the identity on $F$. At each point of $\PP^n$ evaluation gives the full matrix space, in which the locus of rank less than $N$ has codimension $N+1$. Since $N+1>n$, the incidence of matrices dropping rank at some point has dimension smaller than the parameter space of scalar matrices. A general such matrix is therefore surjective everywhere, and so is its pullback.

The finitely many nonempty open subsets intersect in the affine space $\mathcal P_N$. A parameter in their intersection gives a surjection with locally free kernel $K$. At each twist $-imD$, Proposition~\ref{prop-surjhom} gives concentration in degree $i-1$ for both terms, and the induced map in that degree is an isomorphism. The long exact sequence then gives $H^j(Y,K(-imD))=0$ for all $j$ and $1\leq i\leq n$.
\end{proof}

\begin{proof}[Proof of Theorem \ref{thm-main}]
We would like to construct an Ulrich bundle on $X$ with respect to $mH$ for all sufficiently large $m$. By Remark \ref{rem-fixedresidue}, it suffices to fix the residue class $a$ mod $q$ and construct an Ulrich bundle for $(Mq+a)H$ for sufficiently large $M$. By Lemma \ref{lem-passtocover}, it suffices to construct an Ulrich bundle for $f^*(mH)$ on a finite flat cover $f: Y \to X$. Then Theorem \ref{thm-constructbundle} and Proposition \ref{prop-actualbundle} construct the desired bundle on the cyclic cover constructed in \S \ref{sec-cover}. This concludes the proof.      
\end{proof}

\section{Wild representation type}

We first transfer families of Ulrich bundles from projective space. We then prove Theorem \ref{thm-ulrichwild} by varying the map in the kernel presentation \eqref{eq-defK}.

\begin{theorem}\label{thm-wildrep}
    Let $X$ be a connected complex projective scheme of pure dimension $n\ge 2$, and let $H_0$ be a very ample Cartier divisor on $X$. Let $d$ be a positive integer. Suppose that $X$ admits an Ulrich bundle for $H_0$. Then $(X, dH_0)$ is Ulrich-wild if $(n-1)(d-1) \geq 2$.
\end{theorem}

\begin{remark}
The bound $(n-1)(d-1) \geq 2$ is sharp. When $n=2$, we need $d\geq 3$ because every Ulrich bundle on $(\PP^2,\OO_{\PP^2}(2))$ is a direct sum of copies of $T_{\PP^2}$ \cite[Remark~5.3 and Theorem~5.4, p.~4697]{CoskunGenc}. For every $n$, the Ulrich bundles on $(\PP^n,\OO_{\PP^n}(1))$ are precisely the trivial bundles, by \cite[Proposition~2.1(c), p.~548]{ES03} applied to the identity projection.
\end{remark}

\begin{proof}[Proof of Theorem \ref{thm-wildrep}]
    We first show that $(\PP^n,\OO_{\PP^n}(d))$ is Ulrich-wild. By \cite[Lemma~7.1, p.~165]{casnati}, building on \cite[Theorem~A(i) and Corollary~3.1(ii)]{FPL}, it suffices to find a simple Ulrich bundle $E_d$ with $\dim\Ext^1(E_d,E_d)\ge 3$. Let $\lambda_d$ be the partition $((d-1)(n-1),(d-1)(n-2),\ldots,0)$ and set $E_d=(\mathbb S_{\lambda_d}Q)(d-1)$. This irreducible homogeneous bundle is simple by Ramanan's theorem \cite{ramanan}. By \cite[Theorem~5.6, pp.~566--567]{ES03}, $\mathbb S_{\lambda_d}Q$ is acyclic at the twists $-1-jd$ for $0\le j<n$. Thus $E_d(-id)$ is acyclic for $1\le i\le n$, so $E_d$ is Ulrich. Its rank is $d^{\binom{n}{2}}$ by \cite[Corollary~5.7, p.~567]{ES03}.

    We claim that $\dim(\Ext^1(E_d,E_d))\geq 10$ for $(n-1)(d-1)\geq 2$. Under this hypothesis, the endomorphism bundle $\mathcal{E}nd(E_d)=E_d^\vee\otimes E_d$ contains $\mathbb{S}_{2,0,\dots,0,-2}Q$ as a summand, as we show below. In the Borel--Bott--Weil computation \cite[Theorem~1.2, pp.~910--911]{EFW}, adding $\rho=(n,n-1,\ldots,0)$ to $(2,0,\ldots,0,-2,0)$ requires one transposition to give a decreasing sequence. Subtracting $\rho$ gives the highest weight $(2,0,\ldots,0,-1,-1)$ of the unique nonzero cohomology group, in degree one. By the Weyl dimension formula its dimension is
    $$\frac{n(n-1)(n+2)(n+3)}{4}\geq 10.$$

    To see that $\mathcal{E}nd(E_d)$ contains $\mathbb{S}_{2, 0, \dots, 0, -2} Q$, let $V_{\lambda_d}$ be the irreducible $\operatorname{GL}_n$-representation with highest weight $\lambda_d$. We need to show that the representation with highest weight ${2, 0, \dots, 0, -2}$ is a summand of $\End(V_{\lambda_d})$. This is a standard consequence of a highest weight computation, which we briefly recall.

    Let $\theta= e_1 - e_n$ be the highest root of $\mathfrak{s} \mathfrak{l}_n$ and let $e_{\theta}$ be a highest root vector. Let $\kappa:\mathfrak{s} \mathfrak{l}_n \to \End(V_{\lambda_d})$ be the action of $\mathfrak{s} \mathfrak{l}_n$. Consider the operator $\Theta= \kappa(e_\theta)\circ\kappa(e_\theta)$. Since $[e_\alpha, e_\theta]=0$ for every positive root $\alpha$, we have $[\kappa(e_\alpha), \Theta]=0$. For an element $h$ of the Cartan subalgebra, one can compute that $[\kappa(h), \Theta]= 2 \theta(h) \Theta$. Thus, if $\Theta \not=0$, it is a highest weight vector of weight $(2, 0, \cdots, 0, -2)$.  We have $\Theta \not=0$ if $$\langle \lambda_d, \theta^{\vee} \rangle = \lambda_{d,1} - \lambda_{d,n} = (n-1)(d-1) \geq 2.$$ 
    
    Now choose an indecomposable direct summand $V$ of an Ulrich bundle for $H_0$. It is again an Ulrich bundle. For a finite linear projection $\pi:X\to\PP^n$ associated to $H_0$, we have $\pi_*V\cong\OO_{\PP^n}^{\oplus c}$ by \cite[Proposition~2.1(c), p.~548]{ES03}. For each Ulrich bundle $G$ on $(\PP^n,\OO_{\PP^n}(d))$, consider the bundle $W_G=V\otimes\pi^*G$ of \cite[Proposition~5.4, p.~566]{ES03}. The projection formula gives $\pi_*W_G\cong G^{\oplus c}$, so $W_G$ is Ulrich for $dH_0$. Applying Krull--Schmidt on $\PP^n$ \cite[Theorem~3, p.~315]{Atiyah} also shows that $W_G\cong W_{G'}$ implies $G\cong G'$.

    Suppose $G$ is indecomposable. Its endomorphism algebra $R=\End(G)$ is finite-dimensional and local with residue field $\mathbb C$. The residue map is the normalized trace $\operatorname{tr}/\rk(G)$, since every element of the radical is nilpotent. To compute $\End_X(W_G)$, view $\pi_*V$ as a module over $\pi_*\OO_X$. Then $\End_X(W_G)$ consists of the matrices in $\operatorname{Mat}_c(R)$ commuting with the action of $\pi_*\OO_X$ on $G^{\oplus c}$. Each entry of a local action matrix on $G^{\oplus c}$ is the corresponding regular-function entry on $\OO_{\PP^n}^{\oplus c}$ times $\operatorname{id}_G$. Taking normalized traces entrywise in the commutation relations therefore gives an algebra map
\[
\End_X(W_G)\longrightarrow \End_X(V).
\]
It is surjective, with section $u\mapsto u\otimes\operatorname{id}_{\pi^*G}$, and its kernel is contained in the nilpotent ideal $\operatorname{Mat}_c(\operatorname{rad}R)$. Since $V$ is indecomposable, $\End_X(V)$ is local, so $\End_X(W_G)$ is local as well. Thus $W_G$ is indecomposable. Pullback and tensor product carry each family on $\PP^n$ to a family over the same base on $X$, proving the theorem.
\end{proof}

Let $P$ be a very ample divisor on an integral complex projective variety $X$ of dimension $n\ge2$. For each rank $r$, put $h_r=rP^n$. Every Ulrich bundle $B$ of rank $r$ is $0$-regular, with Hilbert polynomial $h_r\binom{t+n}{n}$. Thus $B$ is globally generated, $h^0(B)=h_r$, and $H^i(B)=0$ for $i>0$. Let $Q_r$ be the open subset of the Quot scheme of $\OO_X^{h_r}$ with this Hilbert polynomial where the quotient $B$ is an Ulrich bundle and the map on global sections is an isomorphism. Two points of $Q_r$ give isomorphic bundles precisely when they belong to the same orbit of $G_r=\operatorname{GL}_{h_r}$.

Let $I_s\subset Q_s$ denote the indecomposable locus. When $I_s$ is nonempty, let $d_s$ be the largest dimension of an irreducible locally closed subvariety of $Q_s$ contained in $I_s$ and meeting each $G_s$-orbit in at most one point. This maximum exists since $Q_s$ is of finite type, and its universal quotient supplies the corresponding family.

\begin{lemma}\label{lem-framed-ulrich}
For every $r$ with $Q_r\ne\varnothing$ we have
\[
 \dim Q_r-h_r^2\le r\max_{\substack{1\le s\le r\\ I_s\ne\varnothing}}\frac{d_s+1}{s}.
\]
\end{lemma}

\begin{proof}
The complement of $I_s$ is the union of the images of $G_s\times Q_a\times Q_{s-a}$ under direct sum and change of basis for $1\le a<s$, so $I_s$ is constructible. If $I_s$ is nonempty, choose an irreducible locally closed $S\subset I_s$ of dimension $t=\dim I_s$ and set $c=h_s^2+1$. We claim that $t\le d_s+c$, which is immediate if $t\le c$. Otherwise, embed $S$ in projective space and take a general linear section $T$ of codimension $c$. The pairs of distinct points of $S$ in the same orbit form a constructible set of dimension at most $t+h_s^2$ by the $G_s$-action. Containing two distinct points imposes $2c$ conditions on a linear space, so the locus of such pairs in $T\times T$ has dimension at most $t+h_s^2-2c=\dim T-1$. The closure of its first projection in $T$ has dimension at most $\dim T-1$. Its complement meets each orbit at most once and has dimension $t-c$. An irreducible component of this complement of maximal dimension, with its reduced structure, is therefore an irreducible locally closed subvariety of $Q_s$ of dimension $t-c$ as in the definition of $d_s$. Thus $d_s\ge t-c$.

For a partition $r=r_1+\cdots+r_j$ with all $I_{r_i}$ nonempty, the map from $G_r\times\prod_i I_{r_i}$ to $Q_r$ takes the direct sum and changes its basis of global sections. Its fibers contain free orbits of $\prod_i G_{r_i}$, acting by block changes of basis. Its image therefore has dimension at most
\[
 h_r^2+\sum_i(\dim I_{r_i}-h_{r_i}^2)\ \le\ h_r^2+\sum_i(d_{r_i}+1).
\]
Now $\sum_i(d_{r_i}+1)\le r\max_i\frac{d_{r_i}+1}{r_i}$ since $\sum_i r_i=r$. The finitely many partitions cover $Q_r$, as every Ulrich bundle is a direct sum of indecomposable Ulrich bundles.
\end{proof}

\begin{lemma}\label{lem-moving-kernels}
Let $X$ be an integral complex projective variety of dimension $n\ge2$, let $P$ and $T$ be very ample divisors, and let $A$ be a vector bundle on $X$. Write $\mathcal E=A^\vee\otimes A$. Let $Z$ be the complete intersection of two members of $|T|$, and put
\[
 w=\operatorname{rank}\bigl(H^0(X,\mathcal E(T))\longrightarrow H^0(Z,\mathcal E(T)|_Z)\bigr).
\]
Suppose that $2w>h^0(X,\mathcal E)$ and that for some $N$ a surjection $A^{2N}\to A(T)^N$ has Ulrich kernel for $P$. Then there is a constant $c_P>0$ such that for arbitrarily large ranks $s$ there is a family of dimension at least $c_Ps^2$ of pairwise non-isomorphic indecomposable Ulrich bundles of rank $s$ on $(X,P)$ with fixed determinant.
\end{lemma}

\begin{proof}
Put $b=h^0(X,\mathcal E(T))$ and $e=h^0(X,\mathcal E)$. Block sums give such a surjection for every multiple of the original $N$. For each of these $N$, the maps that are surjective with Ulrich kernel form a nonempty open subset of the $2N^2b$-dimensional space of all maps.

Let $s,t$ cut out $Z$. Consider also the block sum of $N$ copies of $(s,t):A^2\to A(T)$. The cokernel of its dual is $(A^\vee(T)\otimes\mathcal I_Z)^N$ by the Koszul sequence. To specialize in the same parameter space, fix a point outside $Z$ and restrict to the open set of maps surjective at that point. Since $X$ is integral, the dual of every map in this open set is injective as a map of sheaves. The cokernel of the universal dual map is therefore flat over this set by the local criterion for flatness \cite[Tag 00ME]{Stacks}. At a surjection with kernel $K$ this cokernel is $K^\vee$. Upper semicontinuity after tensoring with $A$ now gives, on a nonempty open subset of maps with Ulrich kernel,
\[
 \dim\Hom(K,A)\le N h^0(X,\mathcal E(T)\otimes\mathcal I_Z)=N(b-w).
\]

Set $r=N\rk(A)$ and $h_r=rP^n$. Each of these kernels is globally generated with $h_r$ sections, by Ulrich regularity. Choosing a basis of its sections gives a map from the frame bundle over this open subset to $Q_r$. Over a fixed quotient $\OO_X^{h_r}\twoheadrightarrow B$, a presentation with its frame is determined by an inclusion $B\hookrightarrow A^{2N}$ and an isomorphism of its cokernel with $A(T)^N$. Thus the fiber dimension is at most $2N\dim\Hom(B,A)+N^2 e$. Subtracting this bound from the dimension of the frame bundle gives
\[
 \dim Q_r-h_r^2\ge 2N^2 b-2N^2(b-w)-N^2 e=N^2(2w-e).
\]
Put $\gamma=(2w-e)/\rk(A)^2>0$. Lemma~\ref{lem-framed-ulrich} gives, for each such $r$, an $s\le r$ with $(d_s+1)/s\ge\gamma r$. These $s$ are unbounded as $r\to\infty$, since each $d_s$ is finite, and $d_s\ge\gamma rs-1\ge\gamma s^2-1$.

Let $S\subset I_s$ realize $d_s$. The determinant of the universal family gives a morphism $S\to\operatorname{Pic}(X)$. Since $S$ is irreducible, its image lies in one translate of $\operatorname{Pic}^0(X)$. A fiber over a closed point therefore has a component of dimension at least $d_s-\dim\operatorname{Pic}^0(X)$. Its reduction parametrizes pairwise non-isomorphic bundles with the same determinant. For all sufficiently large $s$ in the sequence this dimension is at least $(\gamma/2)s^2$.
\end{proof}

We need one uniform estimate for the summands $G_\epsilon$ from Theorem \ref{thm-constructbundle}. Fix a residue class $a$ modulo $q$ and its maps $f$ and $\pi$, as in \S\ref{sec-cover}.

\begin{lemma}\label{lem-uniform-homs}
Set $F_\epsilon=A_0\otimes\pi^*G_\epsilon$ and $V_\epsilon=f_*F_\epsilon$. There is a constant $C_0$, depending only on the fixed cover and projections, such that
\[
 \dim\Hom_X(V_\epsilon,V_\eta)\le C_0
\]
for every choice of $\delta\ge u_{\max}+5$ and all $M$ sufficiently large for that choice.
\end{lemma}

\begin{proof}
For a fixed integer $k\ge0$, we first bound $\dim\Hom_{\PP^n}(G_\epsilon,G_\eta(k))$. Write $G_\epsilon=\mathbb S_{\lambda}Q$ and $G_\eta=\mathbb S_{\lambda'}Q$ using dominant integral weights, absorbing any twists into the weights via $\det Q\cong\OO(1)$. By the root formula in \S\ref{sec-numerical} and Theorem \ref{thm-constructbundle}, $|\lambda_j-\lambda'_j|\le2$ for each $j$. In the decomposition of $(\mathbb S_{\lambda}Q)^\vee\otimes\mathbb S_{\lambda'}Q(k)$, only summands $\mathbb S_\mu Q$ with $\mu_n\ge0$ contribute global sections, by Borel--Bott--Weil. Such a $\mu$ is a partition of size $|\lambda'|-|\lambda|+nk\le n(k+2)$, where $|\lambda|$ denotes the sum of the entries. By tensor adjunction its multiplicity equals that of the representation of weight $\lambda'+k\mathbf 1$ in $\mathbb S_\lambda\mathbb C^n\otimes\mathbb S_\mu\mathbb C^n$. Since $\mathbb S_\mu\mathbb C^n$ is a summand of $(\mathbb C^n)^{\otimes|\mu|}$, repeated application of Pieri's rule bounds this multiplicity by $n^{|\mu|}$. There are only finitely many such $\mu$, giving a bound depending only on $k$ and $n$.

For the cyclic cover $f$, the coefficient of degree $\ell-1$ in $f_*\OO_Y=\bigoplus_{j=0}^{\ell-1}L^{-j}$ gives a perfect multiplication pairing with values in $L^{-(\ell-1)}$. Finite duality therefore gives $(f_*F_\epsilon)^\vee\cong f_*(F_\epsilon^\vee\otimes L_Y^{\ell-1})$. The evaluation map
\[
 f^*f_*F_\eta\longrightarrow\bigoplus_{\sigma\in\mu_\ell}\sigma^*F_\eta
\]
is an isomorphism away from the branch divisor and is injective, since its source is locally free and the branch equation is a nonzero divisor. Each $\sigma^*F_\eta$ is isomorphic to $F_\eta$, since $A_0$ is pulled back from $X$ and $\sigma$ acts on the coordinates of $\pi$ by a projective linear transformation. The projection formula now gives
\[
 \dim\Hom_X(V_\epsilon,V_\eta)
 \le \ell\, h^0\!\left(\PP^n,G_\epsilon^\vee\otimes G_\eta(\ell-1)\otimes\pi_*\OO_Y\right).
\]
The vector bundle $\pi_*\OO_Y$ embeds in a fixed direct sum of copies of $\OO_{\PP^n}(k_0)$ for some $k_0\ge0$. Applying the first paragraph with $k=k_0+\ell-1$ proves the lemma.
\end{proof}

\begin{proof}[Proof of Theorem \ref{thm-ulrichwild}]
Let $n=\dim X\ge2$. Fix a residue class $a$ and write $m=qM+a$. With $C_0$ as in Lemma \ref{lem-uniform-homs}, choose $\delta\ge u_{\max}+5$ so large that, for $T=\delta qH$,
\[
 H^1(X,\OO_X(-T))=0,\qquad 2\bigl(h^0(X,\OO_X(T))-2\bigr)>2^n C_0.
\]
This is possible by Lemma \ref{lem-easyconsserre} and ampleness. Keep $\delta$ fixed and take $M$ sufficiently large for the construction in \S\S\ref{sec-numerical}--\ref{sec-cohomology}. Pushing forward the kernel sequence \eqref{eq-defK} along $f$ gives a surjection $A^{2N}\to A(T)^N$ with Ulrich kernel for $mH$, where
\[
 A=f_*F=\bigoplus_\epsilon V_\epsilon^{\oplus\nu_\epsilon}.
\]
Since there are at most $2^n$ summands indexed by $\epsilon$, Lemma \ref{lem-uniform-homs} gives
\[
 h^0(X,A^\vee\otimes A)\le C_0\Bigl(\sum_\epsilon\nu_\epsilon\Bigr)^2
 \le 2^n C_0\sum_\epsilon\nu_\epsilon^2.
\]
Two general sections of $\OO_X(T)$ cut out a complete intersection $Z$ of codimension two, which is zero-dimensional when $n=2$. Their Koszul sequence gives
\[
 0\longrightarrow \OO_X(-T)\longrightarrow \OO_X^{\oplus 2}
 \longrightarrow \mathcal I_Z(T)\longrightarrow 0.
\]
Since $H^0(X,\OO_X(-T))=H^1(X,\OO_X(-T))=0$, we get $h^0(X,\mathcal I_Z(T))=2$. Thus restriction of sections of $\OO_X(T)$ to $Z$ has rank $h^0(X,\OO_X(T))-2$. The scalar endomorphisms of each $V_\epsilon$ split off by the normalized trace. The $\nu_\epsilon^2$ matrix entries in its diagonal block therefore give independent copies of this restriction map. The restriction map for $A^\vee\otimes A(T)$ has rank at least
\[
 \bigl(h^0(X,\OO_X(T))-2\bigr)\sum_\epsilon\nu_\epsilon^2.
\]
Our choice of $\delta$ gives the strict inequality in Lemma~\ref{lem-moving-kernels}. For each sufficiently large $m$ in the chosen residue class, applying the lemma with $P=mH$ gives the asserted families with fixed determinant. There are only finitely many residue classes, so the conclusion holds for every sufficiently large $m$.
\end{proof}


\begin{thebibliography}{CCHMW17}

\bibitem[Ang26a]{Anghel}
C.~Anghel, \emph{Smooth counterexamples to the Eisenbud--Schreyer--Weyman Ulrich existence problem},
preprint (2026), \href{https://arxiv.org/abs/2609.02718v1}{arXiv:2609.02718v1}.

\bibitem[Ang26b]{AnghelACM}
C.~Anghel, \emph{Polarized varieties without arithmetically Cohen--Macaulay bundles and section rings without graded maximal Cohen--Macaulay modules in characteristic zero},
preprint (2026), \href{https://arxiv.org/abs/2609.15589}{arXiv:2609.15589}.

\bibitem[Ang26c]{AnghelSurvey}
C.~Anghel, \emph{The recent counterexamples to the Eisenbud--Schreyer--Weyman Ulrich problem (2003) and to Hochster's small Cohen--Macaulay conjecture (1973): a survey, with a generalization to arbitrary dimension},
preprint (2026), \href{https://arxiv.org/abs/2610.30027}{arXiv:2610.30027}.

\bibitem[Ant21]{Antonelli}
V.~Antonelli, \emph{Characterization of Ulrich bundles on Hirzebruch surfaces},
Rev. Mat. Complut. \textbf{34} (2021), no.~1, 43--74.

\bibitem[AFO17]{AGO}
M.~Aprodu, G.~Farkas, and A.~Ortega, \emph{Minimal resolutions, Chow forms and Ulrich bundles on K3 surfaces},
J. Reine Angew. Math. \textbf{730} (2017), 225--249.

\bibitem[Ati56]{Atiyah}
M.~F.~Atiyah, \emph{On the Krull--Schmidt theorem with application to sheaves},
Bull. Soc. Math. France \textbf{84} (1956), 307--317.

\bibitem[Bea16]{BeauvilleU}
A.~Beauville, \emph{Ulrich bundles on abelian surfaces},
Proc. Amer. Math. Soc. \textbf{144} (2016), no.~11, 4609--4611.

\bibitem[Bea18]{Beauville}
A.~Beauville, \emph{An introduction to Ulrich bundles},
Eur. J. Math. \textbf{4} (2018), no.~1, 26--36.

\bibitem[BN18]{BorisovNuer}
L.~Borisov and H.~Nuer, \emph{Ulrich bundles on Enriques surfaces},
Int. Math. Res. Not. IMRN \textbf{2018} (2018), no.~13, 4171--4189.

\bibitem[BHU87]{BHU}
J.~P.~Brennan, J.~Herzog, and B.~Ulrich, \emph{Maximally generated Cohen--Macaulay modules},
Math. Scand. \textbf{61} (1987), no.~2, 181--203.

\bibitem[CH12]{CHGS}
M.~Casanellas and R.~Hartshorne, \emph{Stable Ulrich bundles},
with an appendix by F.~Gei\ss{} and F.-O.~Schreyer, Internat. J. Math. \textbf{23} (2012), no.~8, Art.~1250083, 50~pp.

\bibitem[Cas17]{CasnatiSurfaces}
G.~Casnati, \emph{Special Ulrich bundles on non-special surfaces with $p_g=q=0$},
Internat. J. Math. \textbf{28} (2017), no.~8, Art.~1750061.

\bibitem[Cas22]{casnati}
G.~Casnati, \emph{Special Ulrich bundles on regular surfaces with non-negative Kodaira dimension},
Manuscripta Math. \textbf{167} (2022), no.~1--2, 151--172.

\bibitem[CG18]{Casnati}
G.~Casnati and F.~Galluzzi, \emph{Stability of rank two Ulrich bundles on projective K3 surfaces},
Math. Scand. \textbf{122} (2018), no.~2, 239--256.

\bibitem[Cos17]{ecoskun}
E.~Coskun, \emph{A survey of Ulrich bundles},
in \emph{Analytic and algebraic geometry} (A.~Aryasomayajula, I.~Biswas, A.~S.~Morye, and A.~J.~Parameswaran, eds.), Springer, Singapore, 2017, pp.~85--106.

\bibitem[CG17]{CoskunGenc}
E.~Coskun and O.~Genc, \emph{Ulrich bundles on Veronese surfaces},
Proc. Amer. Math. Soc. \textbf{145} (2017), no.~11, 4687--4701.

\bibitem[CKM13]{CKM}
E.~Coskun, R.~S.~Kulkarni, and Y.~Mustopa, \emph{The geometry of Ulrich bundles on del Pezzo surfaces},
J. Algebra \textbf{375} (2013), 280--301.

\bibitem[CCHMW17]{CCHMW}
I.~Coskun, L.~Costa, J.~Huizenga, R.~M.~Mir\'o-Roig, and M.~Woolf, \emph{Ulrich Schur bundles on flag varieties},
J. Algebra \textbf{474} (2017), 49--96.

\bibitem[CH20]{CH}
I.~Coskun and J.~Huizenga, \emph{Brill--Noether problems, Ulrich bundles and the cohomology of moduli spaces of sheaves},
Mat. Contemp. \textbf{47} (2020), 21--72.

\bibitem[CHS25]{CHS}
I.~Coskun, J.~Huizenga, and G.~Smith, \emph{Stability and cohomology of kernel bundles on $\PP^n$},
Michigan Math. J. \textbf{75} (2025), no.~1, 173--198.

\bibitem[CNY23]{CNY}
I.~Coskun, H.~Nuer, and K.~Yoshioka, \emph{The cohomology of the general stable sheaf on a K3 surface},
Adv. Math. \textbf{426} (2023), Art.~109102, 85~pp.

\bibitem[CNY25]{CNY2}
I.~Coskun, H.~Nuer, and K.~Yoshioka, \emph{Weak Brill--Noether on abelian surfaces},
Selecta Math. (N.S.) \textbf{31} (2025), no.~3, Art.~48, 31~pp.

\bibitem[CM15]{CostaMiroRoig1}
L.~Costa and R.~M.~Mir\'o-Roig, \emph{$\mathrm{GL}(V)$-invariant Ulrich bundles on Grassmannians},
Math. Ann. \textbf{361} (2015), no.~1--2, 443--457.

\bibitem[CMP12]{CMP}
L.~Costa, R.~M.~Mir\'o-Roig, and J.~Pons-Llopis, \emph{The representation type of Segre varieties},
Adv. Math. \textbf{230} (2012), no.~4--6, 1995--2013.

\bibitem[CLS11]{CLS}
D.~A.~Cox, J.~B.~Little, and H.~K.~Schenck, \emph{Toric varieties},
Graduate Studies in Mathematics, vol.~124, American Mathematical Society, Providence, RI, 2011.

\bibitem[DM17]{DM}
H.~Derksen and V.~Makam, \emph{Polynomial degree bounds for matrix semi-invariants},
Adv. Math. \textbf{310} (2017), 44--63.


\bibitem[EFW11]{EFW}
D.~Eisenbud, G.~Fl{\o}ystad, and J.~Weyman, \emph{The existence of equivariant pure free resolutions},
Ann. Inst. Fourier (Grenoble) \textbf{61} (2011), no.~3, 905--926.

\bibitem[ES03]{ES03}
D.~Eisenbud and F.-O.~Schreyer, \emph{Resultants and Chow forms via exterior syzygies},
with an appendix by J.~Weyman, J. Amer. Math. Soc. \textbf{16} (2003), no.~3, 537--579.

\bibitem[ES09]{ES09}
D.~Eisenbud and F.-O.~Schreyer, \emph{Betti numbers of graded modules and cohomology of vector bundles},
J. Amer. Math. Soc. \textbf{22} (2009), no.~3, 859--888.

\bibitem[ES11]{ES11}
D.~Eisenbud and F.-O.~Schreyer, \emph{Boij--S\"oderberg theory},
in \emph{Combinatorial aspects of commutative algebra and algebraic geometry} (G.~Fl{\o}ystad, T.~Johnsen, and A.~L.~Knutsen, eds.), Abel Symp., vol.~6, Springer, Berlin, 2011, pp.~35--48.


\bibitem[Fae19]{Faenzi2}
D.~Faenzi, \emph{Ulrich bundles on K3 surfaces},
Algebra Number Theory \textbf{13} (2019), no.~6, 1443--1454.

\bibitem[FP21]{FPL}
D.~Faenzi and J.~Pons-Llopis, \emph{The Cohen--Macaulay representation type of projective arithmetically Cohen--Macaulay varieties},
\'Epijournal G\'eom. Alg\'ebrique \textbf{5} (2021), Art.~8, 37~pp.

\bibitem[FR04]{FR}
M.~Fortin and C.~Reutenauer, \emph{Commutative/noncommutative rank of linear matrices and subspaces of matrices of low rank},
S\'em. Lothar. Combin. \textbf{52} (2004), Art.~B52f, 11~pp.

\bibitem[Fuj07]{FujinoToric}
O.~Fujino, \emph{Multiplication maps and vanishing theorems for toric varieties},
Math. Z. \textbf{257} (2007), no.~3, 631--641.

\bibitem[GL26]{GencLopez}
O.~Genc and A.~F.~Lopez, \emph{On the existence of polarized varieties without weakly Ulrich bundles, ordinary instantons with bounded quantum number and supernatural bundles with bounded roots},
preprint (2026), \href{https://arxiv.org/abs/2610.36231}{arXiv:2610.36231}.

\bibitem[Har82]{harris}
J.~Harris, \emph{Curves in projective space},
with the collaboration of D.~Eisenbud, S\'eminaire de Math\'ematiques Sup\'erieures, vol.~85, Presses de l'Universit\'e de Montr\'eal, Montr\'eal, 1982.

\bibitem[Har77]{Hartshorne}
R.~Hartshorne, \emph{Algebraic geometry},
Graduate Texts in Mathematics, vol.~52, Springer-Verlag, New York, 1977.

\bibitem[HUB91]{BaHU}
J.~Herzog, B.~Ulrich, and J.~Backelin, \emph{Linear maximal Cohen--Macaulay modules over strict complete intersections},
J. Pure Appl. Algebra \textbf{71} (1991), no.~2--3, 187--202.

\bibitem[Hor64]{Horrocks}
G.~Horrocks, \emph{Vector bundles on the punctured spectrum of a local ring},
Proc. London Math. Soc. (3) \textbf{14} (1964), no.~4, 689--713.

\bibitem[Hui13]{HuizengaSteiner}
J.~Huizenga, \emph{Restrictions of Steiner bundles and divisors on the Hilbert scheme of points in the plane},
Int. Math. Res. Not. IMRN (2013), no.~21, 4829--4873. \href{https://doi.org/10.1093/imrn/rns195}{doi:10.1093/imrn/rns195}.

\bibitem[Kad21]{Kadets}
B.~Kadets, \emph{Sectional monodromy groups of projective curves},
J. Lond. Math. Soc. (2) \textbf{103} (2021), no.~1, 314--335.

\bibitem[Kim16]{Kim}
Y.~Kim, \emph{Ulrich bundles on rational surfaces with an anticanonical pencil},
Manuscripta Math. \textbf{150} (2016), no.~1--2, 99--110.

\bibitem[Lop26]{Lopez}
A.~F.~Lopez, \emph{On the existence of smooth varieties of any dimension carrying no Ulrich bundles for some very ample line bundle},
preprint (2026), \href{https://arxiv.org/abs/2609.06678v3}{arXiv:2609.06678v3}.





\bibitem[Mil86]{MilneAV}
J.~S.~Milne, \emph{Abelian varieties}, in \emph{Arithmetic geometry} (Storrs, Conn., 1984),
Springer, New York, 1986, pp.~103--150.
\href{https://www.jmilne.org/math/xnotes/AVs.pdf}{Author's corrected version}.

\bibitem[O'G96]{ogrady}
K.~G.~O'Grady, \emph{Moduli of vector bundles on projective surfaces: some basic results},
Invent. Math. \textbf{123} (1996), no.~1, 141--207.

\bibitem[OSS11]{OSS}
C.~Okonek, M.~Schneider, and H.~Spindler, \emph{Vector bundles on complex projective spaces},
with an appendix by S.~I.~Gelfand, Modern Birkh\"auser Classics, Birkh\"auser, Basel, 2011.

\bibitem[Pin26]{Pine}
J.~Pine, \emph{Existence of Ulrich bundles on polarized surfaces and weakly Ulrich bundles on polarized varieties in positive characteristic},
preprint (2026), \href{https://arxiv.org/abs/2610.28456v3}{arXiv:2610.28456v3}.

\bibitem[PS05]{PS}
G.~P.~Pirola and E.~Schlesinger, \emph{Monodromy of projective curves},
J. Algebraic Geom. \textbf{14} (2005), no.~4, 623--642.

\bibitem[Ram66]{ramanan}
S.~Ramanan, \emph{Holomorphic vector bundles on homogeneous spaces},
Topology \textbf{5} (1966), no.~2, 159--177.

\bibitem[Stacks]{Stacks}
The Stacks Project Authors, \emph{The Stacks Project},
2026, \url{https://stacks.math.columbia.edu}.

\end{thebibliography}
\end{document}